
\magnification\magstep1
\baselineskip14.96pt

\newread\AUX\immediate\openin\AUX=\jobname.aux
\newcount\relFnno
\def\ref#1{\expandafter\edef\csname#1\endcsname}
\ifeof\AUX\immediate\write16{\jobname.aux gibt es nicht!}\else
\input \jobname.aux
\fi\immediate\closein\AUX



\def\ignore{\bgroup
\catcode`\;=0\catcode`\^^I=14\catcode`\^^J=14\catcode`\^^M=14
\catcode`\ =14\catcode`\!=14\catcode`\"=14\catcode`\#=14\catcode`\$=14
\catcode`\&=14\catcode`\'=14\catcode`\(=14\catcode`\)=14\catcode`\*=14
\catcode`+=14\catcode`\,=14\catcode`\-=14\catcode`\.=14\catcode`\/=14
\catcode`\0=14\catcode`\1=14\catcode`\2=14\catcode`\3=14\catcode`\4=14
\catcode`\5=14\catcode`\6=14\catcode`\7=14\catcode`\8=14\catcode`\9=14
\catcode`\:=14\catcode`\<=14\catcode`\==14\catcode`\>=14\catcode`\?=14
\catcode`\@=14\catcode`\A=14\catcode`\B=14\catcode`\C=14\catcode`\D=14
\catcode`\E=14\catcode`\F=14\catcode`\G=14\catcode`\H=14\catcode`\I=14
\catcode`\J=14\catcode`\K=14\catcode`\L=14\catcode`\M=14\catcode`\N=14
\catcode`\O=14\catcode`\P=14\catcode`\Q=14\catcode`\R=14\catcode`\S=14
\catcode`\T=14\catcode`\U=14\catcode`\V=14\catcode`\W=14\catcode`\X=14
\catcode`\Y=14\catcode`\Z=14\catcode`\[=14\catcode`\\=14\catcode`\]=14
\catcode`\^=14\catcode`\_=14\catcode`\`=14\catcode`\a=14\catcode`\b=14
\catcode`\c=14\catcode`\d=14\catcode`\e=14\catcode`\f=14\catcode`\g=14
\catcode`\h=14\catcode`\i=14\catcode`\j=14\catcode`\k=14\catcode`\l=14
\catcode`\m=14\catcode`\n=14\catcode`\o=14\catcode`\p=14\catcode`\q=14
\catcode`\r=14\catcode`\s=14\catcode`\t=14\catcode`\u=14\catcode`\v=14
\catcode`\w=14\catcode`\x=14\catcode`\y=14\catcode`\z=14\catcode`\{=14
\catcode`\|=14\catcode`\}=14\catcode`\~=14\catcode`\^^?=14
\Ignoriere}
\def\Ignoriere#1\;{\egroup}

\newcount\itemcount
\def\resetitem{\global\itemcount0}\resetitem
\newcount\Itemcount
\Itemcount0
\newcount\maxItemcount
\maxItemcount=0

\def\FILTER\fam\itfam\tenit#1){#1}

\def\Item#1{\global\advance\itemcount1
\edef\TEXT{{\it\romannumeral\itemcount)}}%
\ifx?#1?\relax\else
\ifnum#1>\maxItemcount\global\maxItemcount=#1\fi
\expandafter\ifx\csname I#1\endcsname\relax\else
\edef\testA{\csname I#1\endcsname}
\expandafter\expandafter\def\expandafter\testB\testA
\edef\testC{\expandafter\FILTER\testB}
\edef\testD{\csname0\testC0\endcsname}\fi
\edef\testE{\csname0\romannumeral\itemcount0\endcsname}
\ifx\testD\testE\relax\else
\immediate\write16{I#1 hat sich geaendert!}\fi
\expandwrite\AUX{\neverexpand\ref{I#1}{\TEXT}}\fi
\item{\ifx?#1?\relax\else\marginnote{I#1}\fi\TEXT}}

\def\today{\number\day.~\ifcase\month\or
  Januar\or Februar\or M{\"a}rz\or April\or Mai\or Juni\or
  Juli\or August\or September\or Oktober\or November\or Dezember\fi
  \space\number\year}
\font\sevenex=cmex7
\scriptfont3=\sevenex
\font\fiveex=cmex10 scaled 500
\scriptscriptfont3=\fiveex
\def\A{{\bf A}}
\def\G{{\bf G}}
\def\P{{\bf P}}

\def\XS{{\widetilde X}}

\def\phi{\varphi}
\def\epsilon{\varepsilon}
\def\theta{\vartheta}
\def\uauf{\lower1.7pt\hbox to 3pt{%
\vbox{\offinterlineskip
\hbox{\vbox to 8.5pt{\leaders\vrule width0.2pt\vfill}%
\kern-.3pt\hbox{\lams\char"76}\kern-0.3pt%
$\raise1pt\hbox{\lams\char"76}$}}\hfil}}

\font\BF=cmbx10 scaled \magstep2

\def\title#1{\par
{\baselineskip1.5\baselineskip\rightskip0pt plus 5truecm
\leavevmode\vskip0truecm\noindent\BF #1\par}
\vskip1truecm
\leftline{\font\CSC=cmcsc10
{\CSC Friedrich Knop}}
\leftline{Emmy-Noether-Zentrum}
\leftline{Department Mathematik, Universit\"at Erlangen, Bismarckstr. 1$1\over2$, D-91054 Erlangen, Germany} 
\leftline{knop@mi.uni-erlangen.de}
\vskip1truecm
\par}

\def\cite#1{\expandafter\ifx\csname#1\endcsname\relax
{\bf?}\immediate\write16{#1 ist nicht definiert!}\else\csname#1\endcsname\fi}
\def\expandwrite#1#2{\edef\next{\write#1{#2}}\next}
\def\neverexpand{\noexpand\noexpand\noexpand}
\def\strip#1\ {}
\def\ncite#1{\expandafter\ifx\csname#1\endcsname\relax
{\bf?}\immediate\write16{#1 ist nicht definiert!}\else
\expandafter\expandafter\expandafter\strip\csname#1\endcsname\fi}
\newwrite\AUX
\immediate\openout\AUX=\jobname.aux
\font\eightrm=cmr8\font\sixrm=cmr6
\font\eighti=cmmi8
\font\eightit=cmti8
\font\eightbf=cmbx8
\font\eightcsc=cmcsc10 scaled 833
\def\eightpoint{%
\textfont0=\eightrm\scriptfont0=\sixrm\def\rm{\fam0\eightrm}%
\textfont1=\eighti
\textfont\bffam=\eightbf\def\bf{\fam\bffam\eightbf}%
\textfont\itfam=\eightit\def\it{\fam\itfam\eightit}%
\def\csc{\eightcsc}%
\setbox\strutbox=\hbox{\vrule height7pt depth2pt width0pt}%
\normalbaselineskip=0,8\normalbaselineskip\normalbaselines\rm}
\newcount\absFnno\absFnno1
\write\AUX{\relFnno1}
\newif\ifMARKE\MARKEtrue
{\catcode`\@=11
\gdef\footnote{\ifMARKE\edef\@sf{\spacefactor\the\spacefactor}\/%
$^{\cite{Fn\the\absFnno}}$\@sf\fi
\MARKEtrue
\insert\footins\bgroup\eightpoint
\interlinepenalty100\let\par=\endgraf
\leftskip=0pt\rightskip=0pt
\splittopskip=10pt plus 1pt minus 1pt \floatingpenalty=20000\smallskip
\item{$^{\cite{Fn\the\absFnno}}$}%
\expandwrite\AUX{\neverexpand\ref{Fn\the\absFnno}{\neverexpand\the\relFnno}}%
\global\advance\absFnno1\write\AUX{\advance\relFnno1}%
\bgroup\strut\aftergroup\@foot\let\next}}
\skip\footins=12pt plus 2pt minus 4pt
\dimen\footins=30pc
\output={\plainoutput\immediate\write\AUX{\relFnno1}}
\newcount\Abschnitt\Abschnitt0
\def\beginsection#1. #2 \par{%
\advance\Abschnitt1
\edef\Aname{\number\Abschnitt}
\vskip0pt plus.10\vsize\penalty-250
\vskip0pt plus-.10\vsize\bigskip\vskip\parskip
\expandafter\ifx\csname#1\endcsname\Aname\relax\else
\immediate\write16{#1 hat sich geaendert!}\fi
\expandwrite\AUX{\neverexpand\ref{#1}{\Aname}}
\leftline{\marginnote{#1}\bf\Aname. \ignorespaces#2}%
\nobreak\smallskip\noindent\SATZ1\GNo0%
\write\TOC{\string\toc{\Aname. #2}{\the\pageno}}}
\def\References{%
\vskip0pt plus.10\vsize\penalty-250
\vskip0pt plus-.10\vsize\bigskip\vskip\parskip
\leftline{\bf References}\nobreak\smallskip\noindent%
\write\TOC{\string\toc{References}{\the\pageno}}}
\newcount\appcount

\def\Appendix#1. #2\par{%
\appcount0
\vskip0pt plus.10\vsize\penalty-250
\vskip0pt plus-.10\vsize\bigskip\bigskip\bigskip\vskip\parskip
\xdef\AppName{#1}
\noindent{\BF Appendix #1: #2}\par
\nobreak\medskip\noindent
\write\TOC{\string\toc{Appendix #1: #2}{\the\pageno}}}

\def\appsection#1. #2 \par{%
\advance\appcount1
\edef\Aname{\AppName\number\appcount}
\vskip0pt plus.10\vsize\penalty-250
\vskip0pt plus-.10\vsize\bigskip\vskip\parskip
\expandwrite\AUX{\neverexpand\ref{#1}{A}}
\leftline{\marginnote{#1}\bf\Aname. \ignorespaces#2}%
\nobreak\smallskip\noindent\SATZ1\GNo0%
\write\TOC{\string\toc{\hskip20pt\Aname. #2}{\the\pageno}}}
\newif\ifmarginalnotes\marginalnotesfalse
\newif\ifmarginalwarnings\marginalwarningstrue

\def\marginnote#1{\ifmarginalnotes\hbox to 0pt{\eightpoint\hss #1\ }\fi}

\def\strutdepth{\dp\strutbox}
\def\Randbem#1#2{\ifmarginalwarnings
{#1}\strut
\setbox0=\vtop{\eightpoint
\rightskip=0pt plus 6mm\hfuzz=3pt\hsize=16mm\noindent\leavevmode#2}%
\vadjust{\kern-\strutdepth
\vtop to \strutdepth{\kern-\ht0
\hbox to \hsize{\kern-16mm\kern-6pt\box0\kern6pt\hfill}\vss}}\fi}

\def\Zitat!{\Randbem{\bf?}{\bf Zitat}}

\newcount\SATZ\SATZ1
\def\proclaim #1. #2\par{\ifdim\lastskip<\medskipamount\removelastskip
\medskip\fi
\noindent{\bf#1.\ }{\it#2\Par}
\ifdim\lastskip<\medskipamount\removelastskip\goodbreak\medskip\fi}
\def\rproclaim #1. #2\par{\ifdim\lastskip<\medskipamount\removelastskip
\medskip\fi
\noindent{\bf#1.\ }{\rm#2\Par}
\ifdim\lastskip<\medskipamount\removelastskip\goodbreak\medskip\fi}
\def\Aussage#1{\expandafter\def\csname#1\endcsname##1.{\resetitem
\ifx?##1?\relax\else
\edef\TEST{#1\penalty10000\ \Aname.\number\SATZ}
\expandafter\ifx\csname##1\endcsname\TEST\relax\else
\immediate\write16{##1 hat sich geaendert!}\fi
\expandwrite\AUX{\neverexpand\ref{##1}{\TEST}}\fi
\proclaim {\marginnote{##1}\Aname.\number\SATZ. #1\global\advance\SATZ1}.}}
\def\rAussage#1{\expandafter\def\csname#1\endcsname##1.{\resetitem
\ifx?##1?\relax\else
\edef\TEST{#1\penalty10000\ \Aname.\number\SATZ}
\expandafter\ifx\csname##1\endcsname\TEST\relax\else
\immediate\write16{##1 hat sich geaendert!}\fi
\expandwrite\AUX{\neverexpand\ref{##1}{\TEST}}\fi
\rproclaim {\marginnote{##1}\Aname.\number\SATZ. #1\global\advance\SATZ1}.}}
\Aussage{Theorem}
\Aussage{Proposition}
\Aussage{Corollary}
\Aussage{Lemma}
\def\Proof:{\par\noindent{\it Proof:}}
\rAussage{Definition}
\def\Remark:{\ifdim\lastskip<\medskipamount\removelastskip\medskip\fi
\noindent{\bf Remark:}}
\def\Remarks:{\ifdim\lastskip<\medskipamount\removelastskip\medskip\fi
\noindent{\bf Remarks:}}
\def\Example:{\ifdim\lastskip<\medskipamount\removelastskip\medskip\fi
\noindent{\bf Example:}}
\def\Examples:{\ifdim\lastskip<\medskipamount\removelastskip\medskip\fi
\noindent{\bf Examples:}}
\font\la=lasy10
\def\strich{\hbox{$\vcenter{\hbox
to 1pt{\leaders\hrule height -0,2pt depth 0,6pt\hfil}}$}}
\def\dashedrightarrow{\hbox{%
\hbox to 0,5cm{\leaders\hbox to 2pt{\hfil\strich\hfil}\hfil}%
\kern-2pt\hbox{\la\char\string"29}}}

\def\Bindestrich{\penalty10000-\hskip0pt}
\let\_=\Bindestrich
\def\.{{\sfcode`.=1000.}}

\def\Par{\par}
\def\:={\mathrel{\raise0,9pt\hbox{.}\kern-2,77779pt
\raise3pt\hbox{.}\kern-2,5pt=}}
\def\=:{\mathrel{=\kern-2,5pt\raise0,9pt\hbox{.}\kern-2,77779pt
\raise3pt\hbox{.}}} 
\def\into{\hookrightarrow}
\def\pfeil{\rightarrow}

\def\Pfeil{\longrightarrow}
\def\pf#1{\buildrel#1\over\rightarrow}
\def\Pf#1{\buildrel#1\over\longrightarrow}

\def\Ugleich{\hbox{$\cup$\kern.5pt\vrule depth -0.5pt}}
\def\|#1|{\mathop{\rm#1}\nolimits}
\def\<{\langle}
\def\>{\rangle}
\let\Times=\times
\def\times{\mathop{\Times}}
\let\Otimes=\otimes
\def\otimes{\mathop{\Otimes}}
\catcode`\@=11
\def\hex#1{\ifcase#1 0\or1\or2\or3\or4\or5\or6\or7\or8\or9\or A\or B\or
C\or D\or E\or F\else\message{Warnung: Setze hex#1=0}0\fi}
\def\fontdef#1:#2,#3,#4.{%
\alloc@8\fam\chardef\sixt@@n\FAM
\ifx!#2!\else\expandafter\font\csname text#1\endcsname=#2
\textfont\the\FAM=\csname text#1\endcsname\fi
\ifx!#3!\else\expandafter\font\csname script#1\endcsname=#3
\scriptfont\the\FAM=\csname script#1\endcsname\fi
\ifx!#4!\else\expandafter\font\csname scriptscript#1\endcsname=#4
\scriptscriptfont\the\FAM=\csname scriptscript#1\endcsname\fi
\expandafter\edef\csname #1\endcsname{\fam\the\FAM\csname text#1\endcsname}
\expandafter\edef\csname hex#1fam\endcsname{\hex\FAM}}
\catcode`\@=12 

\fontdef Ss:cmss10,cmss7,cmss10 scaled 500.
\fontdef Fr:eufm10,eufm7,eufm5.
\def\fa{{\Fr a}}\def\fA{{\Fr A}}

\def\fg{{\Fr g}}\def\fG{{\Fr G}}
\def\fh{{\Fr h}}

\def\fk{{\Fr k}}
\def\fl{{\Fr l}}

\def\fp{{\Fr p}}
\def\fq{{\Fr q}}

\def\fs{{\Fr s}}
\def\ft{{\Fr t}}\def\fT{{\Fr T}}

\fontdef bbb:msbm10,msbm7,msbm5.
\fontdef mbf:cmmib10,cmmib7,.
\fontdef msa:msam10,msam7,msam5.
\def\CC{{\bbb C}}

\def\NN{{\bbb N}}\def\PP{{\bbb P}}
\def\QQ{{\bbb Q}}\def\RR{{\bbb R}}

\def\ZZ{{\bbb Z}}
\def\cA{{\cal A}}\def\cC{{\cal C}}
\def\cE{{\cal E}}\def\cF{{\cal F}}
\def\cI{{\cal I}}\def\cK{{\cal K}}\def\cL{{\cal L}}
\def\cO{{\cal O}}\def\cP{{\cal P}}
\def\cQ{{\cal Q}}

\mathchardef\leer=\string"0\hexbbbfam3F
\mathchardef\subsetneq=\string"3\hexbbbfam24
\mathchardef\semidir=\string"2\hexbbbfam6E
\mathchardef\dirsemi=\string"2\hexbbbfam6F
\mathchardef\haken=\string"2\hexmsafam78
\mathchardef\auf=\string"3\hexmsafam10
\let\OL=\overline
\def\overline#1{{\hskip1pt\OL{\hskip-1pt#1\hskip-.3pt}\hskip.3pt}}
\def\aq{{\overline{a}}}


\def\hq{{\overline{h}}}

\def\Sq{{\overline{S}}}
\def\tq{{\overline{t}}}
\def\uq{{\overline{u}}}

\def\wq{{\overline{w}}}

\def\zq{{\overline{z}}}
%
\newdimen\Parindent
\Parindent=\parindent

\def\hang{\hangindent\Parindent}


\abovedisplayskip 9.0pt plus 3.0pt minus 3.0pt
\belowdisplayskip 9.0pt plus 3.0pt minus 3.0pt
\newdimen\Grenze\Grenze2\Parindent\advance\Grenze1em
\newdimen\Breite
\newbox\DpBox
\def\NewDisplay#1
#2$${\Breite\hsize\advance\Breite-\hangindent
\setbox\DpBox=\hbox{\hskip2\Parindent$\displaystyle{%
\ifx0#1\relax\else\eqno{#1}\fi#2}$}%
\ifnum\predisplaysize<\Grenze\abovedisplayskip\abovedisplayshortskip
\belowdisplayskip\belowdisplayshortskip\fi
\global\futurelet\nexttok\WEITER}
\def\WEITER{\ifx\nexttok\qed\expandafter\leftQEDdisplay
\else\leftdisplay\fi}
\def\leftdisplay{\hskip-\hangindent\leftline{\box\DpBox}$$}
\def\leftQEDdisplay{\hskip-\hangindent
\line{\copy\DpBox\hfill\lower\dp\DpBox\copy\QEDbox}%
\belowdisplayskip0pt$$\bigskip\let\nexttok=}
\everydisplay{\NewDisplay}
\newcount\GNo\GNo=0
\newcount\maxEqNo\maxEqNo=0
\def\eqno#1{%
\global\advance\GNo1
\edef\FTEST{$\fam0(\Aname.\number\GNo)$}
\ifx?#1?\relax\else
\ifnum#1>\maxEqNo\global\maxEqNo=#1\fi%
\expandafter\ifx\csname E#1\endcsname\FTEST\relax\else
\immediate\write16{E#1 hat sich geaendert!}\fi
\expandwrite\AUX{\neverexpand\ref{E#1}{\FTEST}}\fi
\llap{\hbox to 40pt{\ifx?#1?\relax\else\marginnote{E#1}\fi\FTEST\hfill}}}

\catcode`@=11
\def\eqalignno#1{\null\!\!\vcenter{\openup\jot\m@th\ialign{\eqno{##}\hfil
&\strut\hfil$\displaystyle{##}$&$\displaystyle{{}##}$\hfil\crcr#1\crcr}}\,}
\catcode`@=12

\newbox\QEDbox
\newbox\nichts\setbox\nichts=\vbox{}\wd\nichts=2mm\ht\nichts=2mm
\setbox\QEDbox=\hbox{\vrule\vbox{\hrule\copy\nichts\hrule}\vrule}
\def\qed{\leavevmode\unskip\hfil\null\nobreak\hfill\copy\QEDbox\medbreak}
\newdimen\HIindent
\newbox\HIbox
\def\setHI#1{\setbox\HIbox=\hbox{#1}\HIindent=\wd\HIbox}
\def\HI#1{\par\hangindent\HIindent\hangafter=0\noindent\leavevmode
\llap{\hbox to\HIindent{#1\hfil}}\ignorespaces}

\newdimen\maxSpalbr
\newdimen\altSpalbr
\newcount\Zaehler


\newif\ifxxx

{\catcode`/=\active

\gdef\beginrefs{%
\xxxfalse
\catcode`/=\active
\def/{\string/\ifxxx\hskip0pt\fi}
\def\TText##1{{\xxxtrue\tt##1}}
\expandafter\ifx\csname Spaltenbreite\endcsname\relax
\def\Spaltenbreite{1cm}\immediate\write16{Spaltenbreite undefiniert!}\fi
\expandafter\altSpalbr\Spaltenbreite
\maxSpalbr0pt
\gdef\alt{}
\def\\##1\relax{%
\gdef\neu{##1}\ifx\alt\neu\global\advance\Zaehler1\else
\xdef\alt{\neu}\global\Zaehler=1\fi\xdef\SigText{##1\the\Zaehler}}
\def\L|Abk:##1|Sig:##2|Au:##3|Tit:##4|Zs:##5|Bd:##6|S:##7|J:##8|xxx:##9||{%
\def\SigText{##2}\global\setbox0=\hbox{##2\relax}
\edef\TEST{[\SigText]}
\expandafter\ifx\csname##1\endcsname\TEST\relax\else
\immediate\write16{##1 hat sich geaendert!}\fi
\expandwrite\AUX{\neverexpand\ref{##1}{\TEST}}
\setHI{[\SigText]\ }
\ifnum\HIindent>\maxSpalbr\maxSpalbr\HIindent\fi
\ifnum\HIindent<\altSpalbr\HIindent\altSpalbr\fi
\HI{\marginnote{##1}[\SigText]}
\ifx-##3\relax\else{##3}: \fi
\ifx-##4\relax\else{##4}{\sfcode`.=3000.} \fi
\ifx-##5\relax\else{\it ##5\/}\fi
\ifx-##6\relax\else\ {\bf ##6}\fi
\ifx-##8\relax\else\ ({##8})\fi
\ifx-##7\relax\else, {##7}\fi
\ifx-##9\relax\else, \TText{##9}\fi\Par}
\def\B|Abk:##1|Sig:##2|Au:##3|Tit:##4|Reihe:##5|Verlag:##6|Ort:##7|J:##8|xxx:##9||{%
\def\SigText{##2}\global\setbox0=\hbox{##2\relax}
\edef\TEST{[\SigText]}
\expandafter\ifx\csname##1\endcsname\TEST\relax\else
\immediate\write16{##1 hat sich geaendert!}\fi
\expandwrite\AUX{\neverexpand\ref{##1}{\TEST}}
\setHI{[\SigText]\ }
\ifnum\HIindent>\maxSpalbr\maxSpalbr\HIindent\fi
\ifnum\HIindent<\altSpalbr\HIindent\altSpalbr\fi
\HI{\marginnote{##1}[\SigText]}
\ifx-##3\relax\else{##3}: \fi
\ifx-##4\relax\else{##4}{\sfcode`.=3000.} \fi
\ifx-##5\relax\else{(##5)} \fi
\ifx-##7\relax\else{##7:} \fi
\ifx-##6\relax\else{##6}\fi
\ifx-##8\relax\else{ ##8}\fi
\ifx-##9\relax\else, \TText{##9}\fi\Par}
\def\Pr|Abk:##1|Sig:##2|Au:##3|Artikel:##4|Titel:##5|Hgr:##6|Reihe:{%
\def\SigText{##2}\global\setbox0=\hbox{##2\relax}
\edef\TEST{[\SigText]}
\expandafter\ifx\csname##1\endcsname\TEST\relax\else
\immediate\write16{##1 hat sich geaendert!}\fi
\expandwrite\AUX{\neverexpand\ref{##1}{\TEST}}
\setHI{[\SigText]\ }
\ifnum\HIindent>\maxSpalbr\maxSpalbr\HIindent\fi
\ifnum\HIindent<\altSpalbr\HIindent\altSpalbr\fi
\HI{\marginnote{##1}[\SigText]}
\ifx-##3\relax\else{##3}: \fi
\ifx-##4\relax\else{##4}{\sfcode`.=3000.} \fi
\ifx-##5\relax\else{In: \it ##5}. \fi
\ifx-##6\relax\else{(##6)} \fi\PrII}
\def\PrII##1|Bd:##2|Verlag:##3|Ort:##4|S:##5|J:##6|xxx:##7||{%
\ifx-##1\relax\else{##1} \fi
\ifx-##2\relax\else{\bf ##2}, \fi
\ifx-##4\relax\else{##4:} \fi
\ifx-##3\relax\else{##3} \fi
\ifx-##6\relax\else{##6}\fi
\ifx-##5\relax\else{, ##5}\fi
\ifx-##7\relax\else, \TText{##7}\fi\Par}
\bgroup
\baselineskip12pt
\parskip2.5pt plus 1pt
\hyphenation{Hei-del-berg Sprin-ger}
\sfcode`.=1000
\References
}}

\def\endrefs{%
\expandwrite\AUX{\neverexpand\ref{Spaltenbreite}{\the\maxSpalbr}}
\ifnum\maxSpalbr=\altSpalbr\relax\else
\immediate\write16{Spaltenbreite hat sich geaendert!}\fi
\egroup\write16{Letzte Gleichung: E\the\maxEqNo}
\write16{Letzte Aufzaehlung: I\the\maxItemcount}}


\def\message#1{\relax}
\input xy
\xyrequire{cmtip}
\xyrequire{matrix}
\xyrequire{arrow}

\SelectTips{cm}{10}
\def\cxymatrix#1{\vcenter{\xymatrix@=15pt{#1}}}

\def\mathcite#1{\expandafter\ifx\csname#1\endcsname\relax\relax\else
\edef\next{\cite{#1}}\expandafter\dollarstrip\next\fi}
\def\dollarstrip$#1${#1}

\def\rho{\varrho}
\def\phi{\varphi}
\def\hat{\widehat}
\def\fin{^{\rm fin}}
\def\sA{{\Ss A}}
\def\sB{{\Ss B}}
\def\sM{{\Ss M}}
\def\Cinf{\cC^\infty}
\def\ri{{\it i)}}
\def\rii{{\it ii)}}
\def\riii{{\it iii)}}

\def\M{{\Ss M}}
\def\reg{{\rm reg}}
\def\fin{^{\rm fin}}
\def\plim#1{\mathop{\vtop{\offinterlineskip\hbox{\kern1pt\rm lim}\vskip2pt\hbox{$\longleftarrow$}\hbox{$\scriptstyle#1$}}}}

\def\i{\sqrt{-1}}

\Aussage{Conjecture}




\title{Automorphisms of multiplicity free Hamiltonian manifolds}

\newwrite\TOC
\immediate\openout\TOC=\jobname.toc

\def\i{i}

{\baselineskip12pt \noindent {\bf Abstract:} Let $M$ be a multiplicity
free Hamiltonian manifold $M$ for a connected compact Lie group $K$
(not necessarily abelian). Let $\cP$ be the momentum polytope of
$M$. We calculate the automorphism of $M$ as a sheaf over $\cP$ and
show that all higher cohomology groups of this sheaf vanish. From
this, and a recent theorem of Losev, we deduce a conjecture of
Delzant: the momentum polytope and the principal isotropy group
determine $M$ up to isomorphism. Moreover, we give a criterion for
when a polytope and a group are afforded by a multiplicity free
manifold.

}

\bigskip

\beginsection Intro. Introduction

Consider a connected compact Lie group $K$ acting on a connected
Hamiltonian manifold $M$. A measure for the complexity of $M$ is half
the dimension of the symplectic reductions of $M$ and it is natural to
study Hamiltonian manifolds with low complexity first, starting with
the case of complexity zero, the so\_called {\it multiplicity free}
manifolds (see \cite{GS2} or \cite{MiFo}). It has been a longstanding
problem to classify multiplicity free manifolds and it is the purpose
of this paper to complete this project.

More specifically, Delzant conjectured in 1989 that any compact
multiplicity free space is uniquely determined by two invariants: its
momentum polytope $\cP$ and its principal isotropy group
$L_0$. Evidence for this conjecture was Delzant's celebrated
classification of multiplicity free torus actions~\cite{DelTor}, as
well as further particular cases settled by Igl\'esias ($K=SO(3)$,
\cite{Igl}), Delzant ($\|rk|K=2$, \cite{Del2}), and Woodward
(transversal actions, \cite{WoodTrans}). The main objective of this
paper is, building upon work of Losev \cite{Los}, to complete the
proof of Delzant's conjecture (see \cite{DelCon}).

Once we know that a multiplicity free manifold is characterized by the
combinatorial data $(\cP,L_0)$ it is natural to ask which pairs
actually arise this way. In section~\cite{Classification} we show that
this can be reduced to a purely local problem on $\cP$. More
precisely, $\cP$ has to ``look'' locally like the weight monoid of a
smooth affine spherical variety (see \cite{mfClass}). Since the latter
class of varieties has been previously classified by Van Steirteghem
and the author \cite{KVS} this finishes the classification of
multiplicity free manifolds.

The proof of the Delzant Conjecture proceeds in two separate steps:
first a local statement (ultimately due to Losev \cite{Los}) and then a
local-to-global argument (addressed in this paper).

First, we describe briefly the local problem. Let $\ft\subseteq\fk$ be
a Cartan subalgebra. Then it is well known that the orbit space
$\fk^*/K$ can be identified with a Weyl chamber
$\ft^+\subseteq\ft$. Thus the moment map $m:M\pfeil\fk^*$ gives rise to
the {\it invariant moment map} $\psi:M\pfeil\ft^+$. By a celebrated
theorem of Kirwan~\cite{Kir}, the image $\cP=\psi(M)$ is a convex
polytope for $M$ compact. The local statement asserts now that two
compact multiplicity free manifolds with the same momentum polytope
$\cP$ and the same principal isotropy group are isomorphic {\it
locally over $\cP$} (see \cite{localDelzant}). Using techniques from
\cite{Sja}, one can reduce this local problem to a statement about
smooth affine spherical varieties, the ``Knop Conjecture'' (see
\cite{KnopConj}), which was recently settled affirmatively by Losev
\cite{Los}.

To pass from local to global, we need to determine the automorphism
group of a multiplicity free manifold. More precisely, we need to know
the {\it sheaf of automorphisms} $\cA_M$ of $M$ over $\cP$. Our main
result is (see section \cite{Cohomology}):

\Theorem LocGlob. Let $M$ be a compact multiplicity free manifold with
moment polyhedron $\cP$. Then $\cA_M$ is a sheaf of abelian
groups. Moreover, all of its higher cohomology groups vanish:
$H^i(\cP,\cA_M)=0$ for $i\ge1$.

\noindent As an application, the vanishing of $H^1$ implies that two
compact multiplicity free manifolds which are locally isomorphic over
$\cP$ are isomorphic globally. Together with the local statement this
yields the Delzant Conjecture. Moreover, the vanishing of $H^2$
implies that there are no obstructions for gluing manifolds which are
given locally over $\cP$ to one global manifold $M$.

The bulk of this paper is devoted to computing the sheaf $\cA_M$. More
precisely, we show that $\cA_M$ is controlled by a certain root system
$\Phi_M$. In a sense, this root system is a symplectic analogue to the
restricted root system of a symmetric space and is just as fundamental to
understanding the geometry of a multiplicity free manifold.

The main technique of our approach rests on the fact that multiplicity
free manifolds are modeled locally by smooth affine spherical
varieties. This idea is not new, e.g., Sjamaar~\cite{Sja} uses it to
reprove Kirwan's convexity theorem. But instead of (locally) embedding
$M$ into a smooth affine spherical variety $X$ as an {\it open}
subset, we embed $M$ into the cotangent bundle $T^*_X$ as a {\it
closed totally real} subset. This way, one can think of $T^*_X$ as a
complexification of $M$. The geometry of these cotangent bundles has
been the focus of much of our previous research. In particular,
invariants and the automorphisms have been worked out in detail
in~\cite{WuM} and~\cite{ARC}, respectively. Thus, the main task solved
in this paper is to port these algebraic results to the smooth
category.

As a preliminary step for determining all automorphism we have to
first study $K$\_invariant smooth functions on $M$ since these
generate automorphisms via Hamiltonian flows. It will turn out (see
\cite{locWeyl}) that the smooth $K$\_invariants are controlled by a
finite reflection group $W_M$. This group is then used to construct
the root system $\Phi_M$. Note, that this part of the paper is a
specialization of our paper \cite{KnopColl} where manifolds of
arbitrary complexity were considered. We decided not to refer to that
paper, though, since some argument simplify dramatically in the
multiplicity free case.

We conclude this introduction with some historical remarks. The
local-to-global principle \cite{LocGlob} was proved by us around 1995
but circumstances prevented our publishing a formal proof until
now. Around the same time we announced the local statement
\cite{localDelzant} or more precisely its algebraic equivalent
(\cite{KnopConj}) as a conjecture (see e.g. \cite{KnopVerm}). In joint
work with Bart Van Steirteghem, we worked out a
classification~\cite{KVS} of smooth affine spherical varieties in an
attempt to prove the conjecture. Unfortunately, this turned out to be
unfeasible due to the multitude of cases to consider. Meanwhile
(around 2000), Luna launched a program to classify {\it all} spherical
varieties and completed it for groups of type $\sA$~\cite{Lun}. This
enabled his student Camus~\cite{Cam} to settle the ``Knop Conjecture''
for groups of type $\sA$. Finally, Losev~\cite{Los} managed to bypass
all problems which still exist in Luna's program and proved the
conjecture in full generality. This event reinvigorated our interest in
the subject so that finally the proof of Delzant's conjecture is
completely documented.

\medskip\noindent {\bf Acknowledgment:} I would like to thank Yael
Karshon, Eugen Lerman, Reyer Sjamaar, Sue Tolman, and Chris Woodward
for very fruitful discussions. This is especially true for the
analytical and topological side of this paper. At the time of its
conception, I was supported by grants from the NSF and the NSA.

\medskip\noindent {\bf Notation:} In the following, $K$ is a compact
connected Lie group with Lie algebra $\fk$. Let $T_\RR\subseteq K$ be
a maximal torus with Lie algebra $\ft_\RR$ and Weyl group $W$. We fix
moreover a Weyl chamber $\ft^+\subseteq\ft^*_\RR$.

{\it Complexifications:} For any compact Lie group $H$, we denote its
complexification by $H^c$. In particular, we put $G=K^c$, a connected
reductive complex algebraic group. Then $T:=T^c_\RR$ is a maximal
torus of $G$ with Lie algebra $\ft=\ft^c_\RR$. Let
$\Lambda_T:=\|Hom|(T,\CC^\Times)$ the character group of $T$. Then
$T=\|Hom|(\Lambda_T,\CC^\Times)$ with
$\RR$\_structure
$$
\chi(\tq)={\overline{\chi(t)\kern-2pt}\kern2pt}^{-1}.
$$
This way, the group of real points $T(\RR)$ coincides with the compact
torus $T_\RR$. Likewise, $\ft^*=\Lambda\otimes_\ZZ\RR$ carries the real
structure
$$
\overline{\chi\otimes z}:=-\,\chi\otimes\zq
$$
such that $\ft^*_\RR$ is identified with
$$
\{a\in\ft^*\mid\aq=a\}=\Lambda_\ZZ\otimes\i\RR.
$$

\beginsection local. The local Delzant conjecture

Besides developing notation which is used throughout this paper, we
describe in this section how to reduce the local problem mentioned in
the introduction to the purely algebraic statement solved by Losev.

Let $M$ a (possibly non\_compact) connected Hamiltonian $K$\_manifold with
moment map $m:M\pfeil\fk^*$. The {\it momentum image} of $M$ is the
set $\cP=m(M)\cap\ft^+$. A theorem of Kirwan~\cite{Kir}
states that $\cP$ is a convex polyhedron when $M$ is compact. There is
another way to look at it: The restriction of the quotient map
$\fk^*\pfeil\fk^*/K$ to $\ft^+$ is a homeomorphism
$\pi_+:\ft^+\pfeil\fk^*/K$. Now we define the {\it invariant moment
map} as the composed map
$$
\psi:M\pf m\fk^*\auf\fk^*/K\Pf{\pi_+^{-1}}\ft^+.
$$
It is important to keep in mind that $\psi$ is, in general, not
differentiable. The polytope $\cP$ is then the image of $\psi$.

\Definition. A Hamiltonian $K$\_manifold $M$ is {\it multiplicity
free} if it is connected and if $\|dim|M/K=\|dim|\cP$, i.e., if
$\psi:M\pfeil\cP$ has discrete fibers.

If $M$ is compact then (by another theorem of Kirwan, see
\cite{KirBook}) all fibers of $\psi$ are connected. Thus, $M$ is
multiplicity free if and only if $\psi$ is a topological quotient.

In our approach it is necessary to also consider non\_compact
manifolds. More specifically, let $U\subseteq\cP$ be open. Then
$M_U:=\psi^{-1}(U)$ is, in general, a non\_compact multiplicity free
manifold and we need to determine its automorphism group, as well.
For arbitrary non\_compact multiplicity free manifolds certain
pathologies may occur, like $M/K\pfeil\cP$ having non\_connected
fibers or being bijective but not a homeomorphism. Thus we restrict
our attention to {\it convex} manifolds in the sense of
\cite{KnConv}. For multiplicity free manifolds the definition boils
down to:

\Definition. A multiplicity free manifold is {\it convex} if
\Item{}the momentum image $\cP$ is convex and
\Item{}the invariant moment map $\psi:M\pfeil\cP$ is proper.

\noindent Under these conditions, one can show that $M/K\pfeil\cP$ is
a homeomorphism and that $\cP$ is locally polyhedral, i.e., for every
$a\in\cP$ there is a polyhedral cone $C_a\subseteq\ft^*_\RR$ and an open
neighborhood $U$ of $a$ in $\ft^*_\RR$ such that
$$
\cP\cap U=(a+C_a)\cap U.
$$
In particular, $\cP$ is locally closed.

The restriction to convexity is quite mild. Clearly, compact
Hamiltonian manifolds are convex. Moreover, every multiplicity free
manifold is, at least, locally convex. So, our theory applies to some
extent even to the general case. {\it For the remainder of
this section we will assume $M$ to be convex and multiplicity free.}

We proceed by recalling some facts about the principal isotropy
group. Let $\fa^0\subseteq\ft^*_\RR$ be the affine subspace spanned by
$\cP$. The interior of $\cP$ inside $\fa^0$ is called its {\it
relative interior} $\cP^0$. It is open and dense in $\cP$. The
centralizer $L_\RR$ of $\fa^0$ (or, equivalently, $\cP$) is a Levi
subgroup of $K$ containing the maximal torus $T$. The generic
structure of $M$ is then described by the following well\_known

\Lemma GenStr. Let $\Sigma:=m^{-1}(\cP^0)$. Then:
\Item{15} $\Sigma$ is a Hamiltonian $L_\RR$\_manifold with moment map
$m|_\Sigma$.
\Item{16} Let $L_0\subseteq L_\RR$ be the kernel of the action of $L_\RR$ on
$\Sigma$. Then $A_\RR:=L_\RR/L_0$ is a torus acting freely on $\Sigma$.
\Item{17} The map $K\times^{L_\RR}\Sigma\pfeil M$ is an open immersion. Thus,
$L_0$ is a principal isotropy group for $K$ acting on $M$.
\Item{18} The complement $M\setminus K\Sigma$ has codimension
$\ge2$. In particular, $K\Sigma$ is dense and $\Sigma$ is connected.

\Proof: Let $F$ be the smallest face of the Weyl chamber $\ft^+$ which
contains $\cP$ and let $F^0$ be its relative interior. Then
$\cP^1:=\cP\cap F^0$ is open in $\cP$ and contains $\cP^0$. Let
$\Sigma^1:=m^{-1}(\cP^1)$. If $\Sigma$ is replaced by $\Sigma^1$ then
all assertions \cite{I15}--\cite{I18} except for the freeness of the
$A$\_action (which is false on $\Sigma^1$) have been shown in, e.g.,
\cite{LMTW} (even for non\_multiplicity free manifolds). Now observe
that $\Sigma^1$ is a convex multiplicity free $A_\RR$\_manifold with
momentum image $\cP^1$. These actions have been studied by
Delzant~\cite{DelTor} and it follows from his theory that the
$A_\RR$\_action is free over $\cP^0$ and that the points
$x\in\Sigma^1$ with $m(x)\in\partial\cP^1$ have an infinite isotropy
group. Since $\Sigma^1\pfeil\cP^1$ is the quotient by $A_\RR$, we
conclude that the preimage of $\partial\cP^1$ in $\Sigma^1$ has
codimension $\ge2$.\qed

\Remark: It is possible to show that $K\Sigma$ is exactly the union of
all $K$\_orbits of maximal dimension. In particular, it coincides with
the open stratum of $M$.

\medskip

Let $\fa_\RR^*$ be the Lie algebra of $A_\RR$. Then the momentum image
of $\Sigma$ is an open subset in a translate of $\fa_\RR^*$ in
$\ft_\RR^*$. This implies that $\fa_\RR^*$ is in turn the linear
subspace of $\ft_\RR^*$ which is parallel to the affine space
$\fa^0$. From that it follows that the momentum image $\cP$ alone
determines the Lie algebras $\fl_\RR$ and $\fl_0$ of $L_\RR$
and $L_0$, respectively: $\fl_\RR$ is the centralizer of $\cP$ and
$\fl_0$ is the set of $\xi\in\fl_\RR$ with $\chi_1(\xi)=\chi_2(\xi)$
for all $\chi_1,\chi_2\in\cP$.  In other words, the only additional
information gained from $L_0$ is in form of the group
$\Lambda_M:=\|Hom|(A_\RR,\CC^\Times)$ considered as a
lattice in $\i\,\fa_\RR^*\subseteq\fa_\RR^*\otimes_\RR\CC$.

Now we state the local statement mentioned in the introduction. Recall
that $M_U=\psi^{-1}(U)$.

\Theorem localDelzant. For $\nu\in\{1,2\}$ let $M_\nu$ be a convex
multiplicity free Hamiltonian $K$\_manifold with invariant moment map
$\psi_\nu$. Assume $\psi_1(M_1)=\psi_2(M_2)=:\cP$ and
$L_0(M_1)=L_0(M_2)$. Then every $a\in\cP$ has a (convex) open
neighborhood $U$ such that $(M_1)_U\cong(M_2)_U$ as
Hamiltonian $K$\_manifold.

First we reduce to the case $a=0$. Let $\ft^+_a$ be the open star of
$\ft^+$ in $a$, i.e., $\ft^+$ with all faces removed which do not
contain $a$. Let $L_\RR\subseteq K$ be the centralizer of $a$. Then
$\fl^+_a:=\|Ad|L_\RR(\ft^+_a)$ is a connected open subset of
$\fl_\RR^*$. Put $\Sigma_\nu:=m_\nu^{-1}(\fl^+_a)$. The cross\_section
theorem (see e.g. \cite{GSj}~2.4) asserts that
$K\times^{L_\RR}\Sigma_\nu\pfeil M_\nu$ is an open embedding of
Hamiltonian manifolds. It follows that the momentum images of $M_\nu$
and $\Sigma_\nu$ (with respect to $K$ and $L_\RR$, respectively) are
equal in a neighborhood of $a$. Moreover the principal isotropy group
$L_0$ does not change. We may replace therefore $M_\nu$ by
$\Sigma_\nu$. Since $a$ is a character of $\fl_\RR$ we may replace the
moment map $m_{\Sigma_\nu}$ of $\Sigma_\nu$ by its translate
$m_{\Sigma_\nu}-a$. Thus, we may assume from now on that $a=0$.

Let, for the moment, $M$ be any convex Hamiltonian manifold and let
$x\in M$ be a point with $m(x)=0$. Then the symplectic slice theorem
(aka.~equivariant Darboux theorem, see \cite{GSj}~2.3) asserts that a
neighborhood of $Kx$ in $M$ is uniquely determined by two data: the
isotropy group $H_\RR:=K_x$ and the symplectic slice $S=(\fk
x)^\perp/(\fk x)$ which is a symplectic representation of
$H_\RR$. Here, $\fk x$ is the tangent space of the orbit $Kx$ in $x$.

It is well\_known that every symplectic representation of a compact
group carries a compatible unitary structure, i.e., such that the
symplectic form is the imaginary part of the Hermitian form. Thus the
action of $H_\RR$ extends to an action of its complexification
$H$. Now put $X:=G\times^HS$ (recall $G=K^c$) which is a smooth affine
complex algebraic $G$\_variety.

To relate $X$ with $M$ we equip $X$ with the structure of a
Hamiltonian $K$\_manifold: first embed $X$ equivariantly into a finite
dimensional $G$\_module $V$ and choose a $K$\_invariant Hermitian
scalar product on $V$. Then $X$ inherits the structure of a Hamiltonian
$K$-manifold from $V$.

\Lemma LocStruc. The orbits $Kx\subseteq M$ and $K/H_\RR\subseteq
G/H\subseteq X$ have convex open $K$\_invariant neighborhoods which
are isomorphic as Hamiltonian $K$\_manifolds. Moreover, the momentum
images of $M$, $X$ and the open subsets agree in a neighborhood of
$0$.

\Proof: It is easy to see that the orbits $Kx$ and $K/H_\RR$ share the
same local data $(H_\RR,S)$. So the existence of isomorphic open sets
follows from the symplectic slice theorem. Moreover, if $U\subseteq M$
is open and $K$\_invariant then $\psi(U)$ is an open subset of
$\cP$. This implies the second assertion.\qed

The lemma allows to replace $M$ by $X$. Recall that $X$ is {\it
spherical} if a Borel subgroup of $G$ has a dense open orbit. This
is equivalent to the ring of functions $\CC[X]$ being a multiplicity free
$G$\_module. Its structure, as a $G$\_module, is determined be the set
$\Xi_X$ of highest weights (the so\_called {\it weight monoid}). Now,
we have the following comparison theorem due to Brion~\cite{Bri}:

\Theorem Spher. Let $X=G\times^HS$ be as above. Then
\Item{} $X$ is a multiplicity free $K$\_manifold if and only if $X$ is
a spherical $G$\_variety.
\Item{} Let $\cP$ be the momentum image  of $X$ and let $\cQ$ be the
convex cone generated by $\Xi_X$. Then $\cP=\i\cQ$.
\Item{} Let $\Lambda_X$ be the lattice determining the principal
isotropy group $L_0$ of $X$. Then $\Lambda_X=\<\Xi_X\>_\ZZ$.
\Item{19} Conversely, $\Xi_X=\cQ\cap\Lambda_X$.

Now we deduce the local Delzant conjecture from the following
theorem. It was conjectured by the author around 1996 (see,
e.g.,~\cite{KnopVerm}) and proved by Losev~\cite{Los} in 2007.

\Theorem KnopConj. Two smooth affine spherical $G$-varieties with
the same weight monoid are $G$\_equivariantly isomorphic.

The local Delzant conjecture follows: As explained above, we may
assume $a=0$. Likewise, we may replace $M_\nu$ by
$X_\nu=G\times^{H^\nu}S_\nu$ where $H_\RR^\nu\subseteq K$ is a closed
subgroup and $S_\nu$ is a unitary representation of $H_\RR^\nu$. Part
\cite{I19} of \cite{Spher} implies that $X_1$ and $X_2$ have the same
weight monoid. Thus, \cite{KnopConj} provides a $G$\_equivariant
isomorphism $\phi:X_1\pf\sim X_2$.

The variety $X_\nu$ contains a unique closed $G$\_orbit, namely $G/H^\nu$.
Thus $\phi$ induces an
isomorphism $G/H^1\pf\sim G/H^2$. Let
$\phi(eH^1)=g_0H^2$ with $g_0\in G$. Then
$H^2=g_0^{-1}H^1g_0$.

Let $\fk=\fh_\RR^\nu\oplus\fp_\nu$ be a $H_\RR^\nu$\_invariant decomposition. Then
it is well known that the map
$$
K\times^{H_\RR^\nu}\fp_\nu\pfeil G/H^\nu:[k,\xi]\mapsto k\|exp|(\i\,\xi)
$$
is an isomorphism of $K$\_manifolds. It follows that, as a
$K$\_manifold, $G/H^\nu$ has a single minimal $K$\_stratum namely
$K/H_\RR^\nu\times\|exp|(\i\,(\fp_\nu^{H_\RR^\nu}))$. The isomorphism
$\phi$ maps the minimal stratum to a minimal stratum. Thus, we get the
decomposition $g_0=ck$ with $k\in K$ and $H_\RR^2=k^{-1}H_\RR^1k$ and where
$c=\|exp|(\i\,\xi)\in G$ centralizes $H_\RR^2$. Now we change the isomorphism
$\phi$ by the automorphism of $X_2=G\times^{H^2}S_2$ which maps
$[g,s]$ to $[gc^{-1},s]$. Thereby we achieve $c=1$ and therefore
$g_0=k\in K$. By replacing $x_2$ with $kx_2$ we even get
$g_0=k=1$. Then $H^1=H^2$. The isomorphism $S_1\cong S_2$
follows since $S_\nu$ is the complex normal space of
$G/H^\nu\subseteq X_\nu$ in $x_\nu$.\qed

\def\M{{\Ss M}}
\def\sM{\M^c}

\beginsection littleWeyl. The local model and its invariants

Let $M$ be a Hamiltonian manifold and $x\in M$ a point with
$m(x)=0$. The slice theorem asserts that a neighborhood of $Kx$ is
determined by two data: the isotropy group $H_\RR=K_x$ and the
space $S=(\fk x)^\perp/(\fk x)$, considered as a symplectic
representation of $H_\RR$. Let $G$ and $H$ be the complexifications of
$K$ and $H_\RR$, respectively. In the previous section, we worked with
the local model $X=G\times^HS$, but from now on it is more convenient to use
$$19
\M:=K\times^{H_\RR}(\fh_\RR^\perp\oplus S).
$$
Since $K$, $H_\RR$, and $S$ are real algebraic varieties, the same holds
for $\M$. Let $\sM:=\|Spec|\CC[M]$ be the complexification of
$\M$. This variety is related to $X$ in the following way:

\Theorem. Let $\M$, $\sM$, and $X$ be as above. Then
\Item{20} $\sM$ is isomorphic to the cotangent bundle $T^*_X$ of
$X$ as a Hamiltonian $G$\_variety. In particular, the diagram of
moment maps commutes:
$$25
\cxymatrix{\M\ar@{^(->}[r]\ar[d]&\sM\ar@{=}[r]&T^*_X\ar[d]\\
\fk^*\ar@{^(->}[rr]&&\fg^*}
$$
\Item{23} $\M$ is multiplicity free if and only
if $X$ is spherical.

\Proof: Let $\<\cdot,\cdot\>$ denote the $H_\RR$\_invariant Hermitian
scalar product on $S$. Then $S$ is a Hamiltonian manifold with
symplectic form $\omega_S(s_1,s_2)=\|Im|\<s_1,s_2\>$ and moment map
$m_S(s)=(\xi\in\fh_\RR\mapsto{1\over2}\omega_S(\xi s,s))$. Thus also
$$
Z_\RR:=T^*_K\times S=K\times\fk^*\times S
$$
is a symplectic $H_\RR$\_manifold with moment map
$$
m(k,\kappa,s)=-\kappa|_{\fh_\RR}+m_S(s)
$$
After choosing an $H_\RR$\_stable complement $\fp_\RR$ of $\fh_\RR$ in
$\fk$ one can extend the linear form $m_S(s)\in\fh^*_\RR$ to the linear form $\ell_s\in\fk^*$ with $\ell_s(\fp)=0$. Then
$$
K\times\fh_\RR^\perp\times S\Pf\sim m^{-1}(0):(k,\kappa,s)\mapsto
(k,\ell_s+\kappa,s)
$$
After taking the quotient by $H_\RR$ we see that $\M$ is isomorphic to
the symplectic reduction of $Z_\RR$ in $0\in\fk^*$.

Now we complexify. The Hermitian scalar product yields an isomorphism
$$
\Sq\Pf\sim S^*:s\mapsto {1\over2\i}\<\cdot,s\>.
$$
Then $S\into S\otimes_\RR\CC=S\oplus\Sq=S\oplus S^*$ and a short
calculation shows that the restriction of the canonical symplectic
form on $S\oplus S^*=T^*_S$ to $S$ is the standard form
$\omega_S=\|Im|\<\cdot,\cdot\>$. Thus, the embedding of $Z_\RR$ into
its complexification
$$
Z=G\times\fg^*\times S\times S^*=T^*_{G\times S}
$$
is compatible with the Hamiltonian structure. Recall that the
symplectic reduction in $0$ of a cotangent bundle is the cotangent
bundle of the quotient. Thus, the symplectic reduction of $Z$ is
$T^*_X$ showing \cite{I20}. Finally, \cite{I23} is just \cite{Spher}
in conjunction with \cite{LocStruc}.\qed

\Remark: Observe that, even though $T^*_X$ is defined over $\RR$ as a
Hamiltonian variety, neither is $X$ nor the projection $T^*_X\pfeil
X$.

\medskip

The $G$\_invariants on $T^*_X$ have been determined in \cite{WuM}.
We proceed by summarizing the most important facts.

\Theorem Spher0. Let $X$ be a smooth affine spherical variety. Let
$\Lambda_X\subseteq\ft^*$ be the subgroup generated by the weight
monoid and let $\fa^*$ be its $\CC$\_span. Then there is a finite
subgroup $W_X\subseteq GL(\fa^*)$ and a morphism
$q:T^*_X\pfeil\fa^*/W_X$ with:
\Item{}$q$ is the categorical quotient of $T^*_X$ by $G$.
\Item{}$W_X$ is generated by reflections. In particular, $\fa^*/W_X$
is isomorphic to an affine space.
\Item{}$W_X$ is a subgroup of
$N_W(\fa^*)/C_W(\fa^*)$ and normalizes the
lattice $\Lambda_X\subseteq\fa^*$.

\Proof: This is \cite{WuM}~Satz~6.4, Satz~6.6 and \cite{CM} Lemma~3.4,
Theorem~4.2.\qed

\noindent
To apply this theorem to $\M$, recall that
$\fa^*=\Lambda_X\otimes_\ZZ\CC$ has the real structure
$$
\overline{\chi\otimes z}=-\chi\otimes\zq.
$$
Thus, the real points are $\fa_\RR^*=\Lambda_X\otimes\i\,\RR$. Recall
(\cite{Spher}) that $\psi:\M\pfeil\ft^+$ takes values in $\fa_\RR^*$.

For any $n\in N_W(\fa^*)$ put ${}^nW_X:=nW_Xn^{-1}$. Then there is
an isomorphism
$$
\tilde n:\fa^*/W_X\Pf\sim\fa^*/{}^nW_X:W_Xa\mapsto{}^nW_Xna.
$$

\Theorem commute. Assume that $\M$ is multiplicity free. Then there is $n\in
N_W(\fa^*)$ such that the following diagram commutes:
$$21
\cxymatrix{
\M\ar[r]\ar[d]_\psi&\sM\ar[d]_{\tilde nq}\\
\fa_\RR^*\ar[r]&\fa^*/W_\M}
$$
Here $W_\M:={}^nW_X$.

\Proof: Consider the fiber product
$$
\sM_N:=\sM\times_{\fa^*/N}\fa^*.
$$
where $N=N_W(\fa^*)/C_W(\fa^*)$. This variety is in general not
irreducible and one of its irreducible components is
$$
\sM_e:=\sM\times_{\fa^*/W_X}\fa^*
$$
(see \cite{CM} Lemma~3.4). The others are translates of $\sM_e$ by an
element $n\in N$:
$$
\sM_n:=\sM\times_{\fa^*/W_X}\kern-10pt{}^n\kern6pt\fa^*
$$
where the superscript ${}^n$ indicates that the map
$\fa^*\pfeil\fa^*/W_X$ is twisted by $\fa^*\pf{n\cdot}\fa^*$.

The morphism $\fa^*/N\pfeil\ft^*/W$ is finite and birational
onto its image. Thus $\fa^*/N$ is the normalization of the image
of $\fa^*$ in $\ft^*/W$. The compatibility of moment maps~\cite{E25}
yields therefore the commutative diagram
$$
\cxymatrix{
&\M\ar[dl]_\psi\ar[dd]\ar@{^(->}[r]&\sM\ar[dd]\\
\fa_\RR^*\ar[dr]\\
&\fa_\RR^*/N\ar@{^(->}[r]&\fa^*/N}
$$
hence a map
$$17
\phi:\M
\pfeil
\M\times_{\fa_\RR^*/N}\fa_\RR^*
\into
\sM\times_{\fa^*/N}\fa^*
$$
Because $\psi$ is not smooth, we cannot expect $\phi$ to be smooth.
But if we restrict to the dense open subset $\M^0:=K\Sigma=Km^{-1}(\cP^0)$ as in
\cite{GenStr} then $\psi$ and $\phi$ become even real analytic. 
Since $\M^0$ is connected, we conclude that $\phi$ maps
$\M^0$ into a single irreducible component $\sM_n$ for some $n\in N$. By
continuity, we obtain $\phi(\M)\subseteq\sM_n$. This means that the diagram
$$
\cxymatrix{\M\ar[d]_\psi\ar@{^(->}[r]&\sM\ar[d]_q\\\fa_\RR^*\ar[r]^-n&\fa^*/W_X}
$$
commutes. This is easily seen to be equivalent to the commutativity of
\cite{E21}.\qed

\Remark: The appearance of $n$ is unavoidable at this stage since we
didn't specify how $q$ was defined. It is very likely, however, that
$n$ is always $1$. With the methods of \cite{Wood} it would be
possible to show that $n$ normalizes $W_X$, i.e., that $W_\M=W_X$.

\medskip

Now we use the preceding theory to determine the smooth
$K$\_invariants on $\M$.

\Corollary surAlg. The composed map
$\psi/W_\M:\M\pf\psi\fa_\RR^*\pfeil\fa_\RR^*/W_\M$ is smooth and
$$26
(\psi/W_\M)^*:\cC^\infty(\fa_\RR^*)^{W_\M}\auf\cC^\infty(\M)^K.
$$
is surjective.

\Proof: The map $\M\pfeil\fa_\RR^*/W_\M$ is a quotient by $K$ in the algebraic
category. A theorem of Schwarz~\cite{Sch} asserts that then
$$
\cC^\infty(\fa_\RR^*/W_\M)\auf\cC^\infty(\M)^K
$$
is surjective. The same theorem applied to $W_\M$ shows that
$\cC^\infty(\fa_\RR^*/W_\M)=\cC^\infty(\fa_\RR^*)^{W_\M}$.\qed

Since $\M$ is just a local model we are going to need a local
refinement. Let $U\subseteq\fa_\RR^*$ be a subset which we do {\it not}
assume to be $W_\M$\_invariant. Then, in abuse of language, we denote
the image of $U$ in $\fa_\RR^*/W_\M$ by $U/W_\M$. It carries an induced
differentiable structure: a function on $U/W_\M$ is smooth if it is
locally the restriction of a smooth function on $\fa_\RR^*/W_\M$.

\Corollary surAlgLoc. For an open subset $U$ of $\cP$, the momentum image of
$\M$, let $\M_U:=\psi^{-1}(U)\subseteq\M$. Then
$$
\cC^\infty(U/W_\M)\Pf\sim\cC^\infty(\M_U)^K
$$

\Proof: Injectivity is clear. To prove surjectivity, let $f$ be some
$K$\_invariant on $\M$. A partition of unity argument shows that we
may assume that $f$ has compact support. But then it extends to all of
$\M$ by zero and we conclude with \cite{surAlg}.\qed

\beginsection infinitesimalAut. Invariants of multiplicity free manifolds

In this section, we determine the smooth invariants on an arbitrary
convex multiplicity free manifold.

\Theorem locWeyl. Let $M$ be a convex multiplicity free Hamiltonian
manifold with momentum image $\cP=\psi(M)$. Let $\fa^0\subseteq\ft_\RR^*$
be the affine space spanned by $\cP$. Then: 
\Item{28} There is a finite group $W_0\subseteq N_W(\fa^0)/C_W(\fa^0)$
such that the composed map $\psi/W_0:M\pfeil\fa^0/W_0$ is smooth and
the induced map
$$30
(\psi/W_0)^*:\cC^\infty(\cP/W_0)\Pf\sim\cC^\infty(M)^K.
$$
is an isomorphism.
\Item{29} Among all groups $W_0$ as in \cite{I28}, there is a unique
minimal one, denoted by $W_M$. It is characterized by the fact that it
is generated by reflections which have a fixed point in $\cP$.

\Remarks: 1. It is easy to construct an example where $\ft_\RR^*=\RR$,
$\ft^+=\RR_{\ge0}$, $W=\{\pm1\}$ and $\cP\subseteq\ft^+$ is a closed
interval not containing the origin (take e.g. $K=SU(2)$ acting on
$\PP^2(\CC)$ blown up at one point). Then the map
$\cP\pfeil\RR/\{\pm1\}$ is a diffeomorphism onto its image. Thus both
$W_0=1$ and $W_0=\{\pm1\}$ would work. In particular, the group $W_0$
as in \cite{I28} is not unique.

2. Let $\tilde W\subseteq N_W(\fa^0)/C_W(\fa^0)$ be the subgroup
generated by all reflections which have a fixed point in $\cP$. Then
$W_M\subseteq\tilde W$ with equality being the rule rather than the
exception. The standard example for when $W_M\subsetneq\tilde W$ is
$K=SU(2)$ acting on $M=\CC^2$. In fact, in that case the moment map is
homogeneous of degree $2$ and the generating invariant on $\fk^*$ is
of degree $2$ as well. So, the minimal degree of an invariant which is
pull\_back from $\ft^*/W$ is $4$. Thus, the degree $2$ invariant
$q(z)=|\!|z|\!|^2$ on $M$ cannot be a pull\_back. This shows $W_M=1$
whereas $\tilde W=W=\ZZ/2\ZZ$.

3. The theorem generalizes to to arbitrary convex Hamiltonian
$K$\_manifolds if one replaces the ring of $K$\_invariants in
\cite{E30} by its Poisson center. See \cite{KnopColl} for details.

\medskip

First, we prove a local version of the theorem.

\Lemma locInv. For every $a\in\cP$ there is a unique subgroup $W(a)$
of $N_W(\fa^0)/C_W(\fa^0)$ having $a$ as a fixed point and an open
convex neighborhood $U$ of $a$ in $\cP$ such that
$$27
\cC^\infty(U/W(a))\Pf\sim\cC^\infty(M_U)^K
$$
where $M_U=\psi^{-1}(U)$. Moreover, $W(a)$ is generated by reflections.

\Proof: We first establish the existence of $W(a)$ and $U$.
Let $Kx=\psi^{-1}(a)$. If $a=0$ then the assertion follows
from \cite{surAlgLoc} and the fact that $M$ is isomorphic to some
model $\M$ near $Kx$. In this case $W(a)=W_\M$.

If $a\ne0$ but $a\in(\fk^*)^K$ then we can shift the moment map by
$-a$. Since $a$ is also fixed by $W$ it follows that
$N_W(\fa^*)/C_W(\fa^*)=N_W(\fa^0)/C_W(\fa^0)$. Thus, we are reduced to
the case $a=0$.

In general, we apply the cross\_section theorem. It states that a
neighborhood $M_0$ of $Kx$ in $M$ is isomorphic to $K\times^L\Sigma$
where $L=K_a$ and $\Sigma$ is a multiplicity free
$L$\_manifold. Moreover, the momentum images of $M_0$ and $\Sigma$
coincide and form an open neighborhood of $a$ in $\cP$. Because of
$\cC^\infty(K\times^L\Sigma)^K=\cC^\infty(\Sigma)^L$ we may replace
$M$ by $\Sigma$. Thus we are reduced to the case $a\in(\fk^*)^K$.

For the uniqueness part, let $W_1$ and $W_2$ be two groups having all
advertised properties. Then there is a neighborhood $U$ of $a$ in
$\cP$ with
$$28
\cC^\infty(U/W_1)=\cC^\infty(M_U)^K=\cC^\infty(U/W_2).
$$
Let $\RR[\![\fa^0]\!]$ be the ring of formal power series in $a$. Then
\cite{E28} implies
$$
\RR[\![\fa^0]\!]^{W_1}=\RR[\![\fa^0]\!]^{W_2}
$$
(here we use $W_\nu a=a$ for this to make sense). But this implies
$W_1=W_2$ (look at fields of fractions and apply Galois theory).\qed

To glue the local Weyl groups $W(a)$ to a global one we need the
following coherence property:

\Lemma Coherence. Every $a\in\cP$ has a neighborhood $U\subseteq\cP$ with
$$5
W(b)=W(a)_b\quad\hbox{for all }b\in U.
$$

\Proof: Note that $\fa^0/W(a)_b\pfeil\fa^0/W(a)$ is an isomorphism
near $b$. Thus $W(a)_b=W(b)$ by uniqueness of $W(b)$.\qed

Next we prove a criterion for gluing all local groups $W(a)$ together:

\Lemma RootCrit. Let $\Delta$ be a finite root system and
$\Delta^+\subseteq\Delta$ a system of positive roots. Let $\Sigma$ be
a set of positive roots such that $\<\alpha\mid\beta^\vee\>\le0$ for
all $\alpha,\beta\in\Sigma$, $\alpha\ne\beta$. Then $\Sigma$ is the
set of simple roots for a subroot system of $\Delta$.

\Proof: Let $C$ be the $|\Sigma|\times|\Sigma|$\_matrix whose entries
are $\<\alpha|\beta^\vee\>$ with $\alpha,\beta\in\Sigma$. Then $C=C'D$
where $C'$ is the Gram matrix $(\alpha,\beta)$ of $\Sigma$, hence
positive semidefinite, and where $D$ is a diagonal matrix with the
positive entries $2/(\beta,\beta)$ on the diagonal. If $C$ were
singular, then \cite{Kac} Thm.~4.3 implies that there is a nonzero
non\_negative vector $v$ such that $Cv=0$. In other words, $0$ would
be a positive linear combination of $\Sigma$ which is impossible since
$\sigma\subseteq\Delta^+$. It follows that $C$ is a generalized Cartan
matrix of finite type. Such a matrix is associated to a finite root
system in a unique way. In particular, $\Sigma$ is its set of simple
roots.\qed

\Lemma globWeyl. There is a subgroup $W_M$ of $N:=N_W(\fa^0)/C_W(\fa^0)$ with
$W(a)=(W_M)_a$ for all $a\in\cP$.

\Proof: We claim that the group $W_M$ which is generated by
$\bigcup_{a\in\cP}W(a)$ has all required properties. From
$W_M\subseteq N$ it follows that $W_M$ finite. Since each individual
$W(a)$ is generated by reflections, $W$ is also generated by
reflections. From $\cP\subseteq\ft^+$ it follows that $W_x=C_W(\fa^0)$
for all $x\in\cP^0$. This means that every reflection hyperplane of
$W_M$ meets $\cP$ at most in its boundary. Hence, $\cP$ is entirely
contained in a single Weyl chamber, say $\fa^+$, of $W_M$.

The finiteness of $W_M$ also implies that it has a fixed point
$a_0\in\fa^0$. Thus we may identify $\fa^0$ with $\fa_\RR^*$ and assume
that $W_M$ is acting linearly.

Since $N$ is a subquotient of $W$ and since $\fa_\RR^*$ is defined over
$\QQ$ there is a lattice in $\fa_\RR^*$ which is normalized by $W_M$. This
means that $W_M$ is in fact a crystallographic reflection group, i.e.,
that it is the Weyl group attached to a root system $\Delta$.

Every $W(a)$ is the Weyl group of a subroot system
$\Delta(a)\subseteq\Delta$.  Moreover, the Weyl chamber $\fa^+$
determines a Weyl chamber for $\Delta(a)$ and therefore a set
$\Sigma(a)\subseteq\Delta(a)$ of simple roots. Moreover, the coherence
property \cite{Coherence} implies
$$20
\Sigma(b)=\{\gamma\in\Sigma(a)\mid\gamma(b)=0\}
$$
for all $a\in\cP$ and all $b\in\cP$ sufficiently close to $a$. This
follows from the fact that the isotropy group $W(a)_b$ is generated
by simple reflections if $b$ is in the dominant Weyl chamber.

Now let $\Sigma$ be the union of all $\Sigma(a)$, $a\in\cP$. We
claim that $\Sigma(a)$ can recovered from $\Sigma$, namely:
$$29
\Sigma(a)=\{\alpha\in\Sigma\mid\alpha(a)=0\}.
$$
Indeed, \cite{E20} implies ``$\subseteq$''. Conversely, let
$\alpha\in\Sigma$ with $\alpha(a)=0$. By construction,
$\alpha\in\Sigma(b)$ for some $b\in\cP$. Thus
$\alpha(a)=\alpha(b)=0$ and therefore $\alpha(c)=0$ for all $c$ on the
line segment $[a,b]$. The convexity of $\cP$ implies
$[a,b]\subseteq\cP$. From \cite{E20} it follows that the set of
$c\in[a,b]$ with $\alpha\in\Sigma(c)$ is both open and closed. Hence
$\alpha\in\Sigma(a)$, proving the claim.

Next, we want to apply \cite{RootCrit}. For that let
$\alpha\in\Sigma(a)$ and $\beta\in\Sigma(b)$ with $a,b\in\cP$. Then
$a,b\in\fa^+$ implies $\alpha(b)\ge0$, $\beta(a)\ge0$. If both
$\alpha(b)>0$ and $\beta(a)>0$ then $\gamma(a)<0$ and $\gamma(b)>0$
where $\gamma:=\alpha-\beta$. Thus $\gamma$ is not a root which
implies $\<\alpha,\beta^\vee\>\le0$. On the other hand if
$\alpha(b)=0$ then $\alpha,\beta\in\Sigma(b)$ by \cite{E29}. Thus,
again $\<\alpha,\beta^\vee\>\le0$. The same conclusion holds for
$\beta(a)=0$. \cite{RootCrit} implies that $\Sigma$ is the set of
simple roots for some subroot system of $\Delta$ and therefore of
$\Delta$ itself since $W_M$ is generated by the reflections in
$\Sigma$. It follows that $(W_M)_a$ is generated by the
reflections corresponding to the roots $\alpha\in\Sigma$ with
$\alpha(a)=0$. Thus $(W_M)_a=W(a)$ by \cite{E29}.\qed

\noindent{\it Proof of \cite{locWeyl}:} Put $W_0=W_M$ as in
\cite{globWeyl}. Then $M\pfeil\fa_\RR^*/W_M$ is smooth since it is locally
smooth. Moreover, by a partition of unity argument, \cite{locInv} shows
that the map \cite{E30} is surjective. Finally, $W_M$ is the minimal
possibility for $W_0$ since $W_0$ has to contain all local Weyl groups
$W(a)$.\qed

We conclude this section with a general statement about smooth
functions on $\cP/W_M$ which might be useful for concrete
calculations. More precisely, since $\cP\pfeil\cP/W_M$ is bijective,
functions on $\cP$ are the same as functions on $\cP/W_M$. So it is
interesting to characterize those functions on $\cP$ which are smooth
on $\cP/W_M$.

\Proposition smoothFunctions. Let $f:\cP\pfeil\RR$ be a function. Then the
following properties are equivalent:
\Item{40}$f\in\cC^\infty(\cP/W_M)$.
\Item{41}For every $a\in\cP$ there is a $(W_M)_a$\_invariant smooth
function $\tilde f$ on $\fa^*$ such that $f$ and $\tilde f$ coincide
in a neighborhood of $a$ in $\cP$.
\Item{42}For every $a\in\cP$ the Taylor series of $f$
in $a$ is $(W_M)_a$\_invariant.\Par\noindent
If $\cP$ is closed in $\fa^*$ then these conditions are equivalent to
\Item{43}$f$ is the restriction of a smooth $W$\_invariant function on
$\fa^*$.
\Item{44}Let $g_1,\ldots,g_r$ be generators of $\RR[\fa^*]^{W_M}$. Then
there is a smooth function $h(x_1,\ldots,x_r)$ with
$f=h(g_1,\ldots,g_r)$ on $\cP$.

\Proof: The implications
\cite{I40}$\Rightarrow$\cite{I41}$\Rightarrow$\cite{I42} are
clear. The implication \cite{I42}$\Rightarrow$\cite{I40} is a
corollary of a deep theorem of Bierstone\_Milman (\cite{BM}
Thm.~3.2). In the case $\cP$ is closed, the last two conditions clearly
imply the first three. Moreover, \cite{I44}$\Rightarrow$\cite{I43}. For
the implication \cite{I40}$\Rightarrow$\cite{I44} observe that
$(g_1,\ldots,g_r)$ embeds $\fa^*/W_M$, hence also $\cP/W_M$, into
$\RR^r$ as a closed semi\_algebraic subset. Then e.g. a partition of
unity argument shows that $f$ can be extended to a smooth function $h$
on $\RR^r$.\qed

\def\fq{\fs}

\beginsection groupscheme. The group scheme

The automorphisms of the cotangent bundle $T^*_X$ of a smooth
spherical variety have been studied in \cite{ARC}. There it was shown
that they can be interpreted as the global sections of a certain
abelian group scheme acting on $T^*_X$. Our program is to extend this
theory to multiplicity free manifolds. 

In the present section, we recall the definition and some of the
fundamental properties of the group scheme. For this, let $\fa$ be any
finite dimensional $\CC$-vector space and let $W\subseteq GL(\fa)$ be
a finite reflection group. We assume that $W$ normalizes a lattice
$\Lambda\subseteq\fa^*$. This implies that $W$ acts on the torus
$A:=\|Hom|(\Lambda,\CC^\Times)$. We have $X(A)=\Lambda$ and
$\|Lie|A=\fa$.  Let $\fq:=\fa^*/W=\|Spec|S^*(\fa)^W$ be the quotient
with quotient morphism $\pi:\fa^*\pfeil\fq$. It is well-known that
$\fq$ is isomorphic to an affine space and that $\pi$ is finite and
faithfully flat.

For any $\fq$-scheme $Y$ let $Y_{\fa^*}$ be the pull back-scheme
$Y\times_\fq\fa^*$. It carries a natural action of $W$ and we
are interested in $W$-equivariant morphisms $Y_{\fa^*}\pfeil A$.

\Proposition. {\rm(\cite{ARC} Lemma~2.3)} The functor $Y\mapsto
\|Mor|^W(Y_{\fa^*},A)$ is represented by a smooth affine abelian group
scheme $\fA_W/\fq$. This means that a $W$-equivariant morphism
$Y_{\fa^*}\pfeil A$ is the ``same'' as a section of $\fA_W$ over $Y$,
i.e., an $\fq$\_morphism $Y\pfeil\fA_W$.

\Example: Let $\fa=\CC$, $W=\{\pm1\}$, and $\Lambda=\ZZ$. If $t$ is
the coordinate function on $\fa^*\cong\CC$ then $s=t^2$ is the
coordinate function on $\fq\cong\CC$. Let $\fA/\fq$ be the subgroup
scheme of $GL(2,\CC)\times\fq$ consisting of all pairs $(M,s)$ with
$$33
M:=\pmatrix{a&sb\cr b&a\cr}
$$
and $\|det|(M)=a^2-sb^2=1$. We claim $\fA_W \cong\fA$. For this, we
verify that $\fA$ has the required universal property. A morphism
$Y\pfeil\fq$ is given the image of $s$ in $\cO(Y)$. Then every
function on $Y_{\fa^*}$ has the form $a+tb$ with $a,b\in\cO(Y)$ and
$t^2=s$. A $W$\_equivariant morphism $\phi:Y_{\fa^*}\pfeil
A=\CC^\Times$ amounts to an invertible function $a+tb$ with
$(a+tb)^{-1}=a-tb$. Thus, $\phi$ is equivalent to a pair
$a,b\in\cO(Y)$ with $a^2-sb^2=(a+tb)(a-tb)=1$, i.e., to a
$\fq$\_morphism $Y\pfeil\fA$.

\medskip

In the sequel, we need a refinement of $\fA_W$. In the example observe
that the generic fiber of $\fA$ is $\fA_s\cong\CC^\Times$ while the
zero fiber is $\fA_0\cong\{\pm1\}\times\CC$. Thus, we obtain an open
subgroup scheme $\fA'\subset\fA$ by removing the set $\{-1\}\times\CC$
inside $\fA_0$ from $\fA$. Since that set is a divisor, the group
scheme $\fA'$ is still affine.

In general, it was shown in \cite{ARC} that the {\it affine} open
subgroup schemes of $\fA_W $ correspond to the root systems with Weyl
group $W$ and lattice $\Lambda$. More precisely, let
$\Phi=(\Lambda,\Delta,\Lambda^\vee,\Delta^\vee)$ be a reduced root
datum. Let $W=W(\Phi)$ be the Weyl group of $\Phi$ and let
$A=\|Hom|(\Lambda,\CC^\Times)$ be the torus with character group
$\Lambda$. Then $\fa^*=\Lambda\otimes_\ZZ\CC$.

For a root $\alpha\in\Delta$ let $s_\alpha\in W$ be the corresponding
reflection. If we regard $\alpha$ as a character $A\pfeil\CC^\Times$ and the
coroot $\alpha^\vee$ as a cocharacter $\CC^\Times\pfeil A$ then the action of $s_\alpha$ on $A$ is given by the formula
$$
s_\alpha(t)=t\cdot\alpha^\vee(\alpha(t))^{-1}
$$
Let $\chi\in\Lambda$ be an arbitrary character and let $t\in A$ with
$s_\alpha(t)=t$. Then
$$14
1=\chi\big(t\cdot s_\alpha(t)^{-1}\big)=\alpha(t)^{\<\chi\mid \alpha^\vee\>}
$$
The choice $\chi=\alpha$ implies in particular $\alpha(t)\in\{1,-1\}$.

Now let $Z$ be a $W$-scheme and $Z^{s_\alpha}$ the $s_\alpha$\_fixed
point scheme. Let $\phi:Z\pfeil A$ be a $W$\_equivariant
morphism. Then the discussion shows
$\alpha(\phi(Z^{s_\alpha}))\subseteq\{1,-1\}$.

\Definition. A morphism $\phi:Z\pfeil A$ is called
{\it $\Phi$\_equivariant} if it is $W$\_equivariant and if
$$15
\alpha\big(\phi(Z^{s_\alpha})\big)=1
$$
for all roots $\alpha\in\Delta$.

\Remark: If there exists $\chi\in\Lambda$ with $\<\chi\mid
\alpha^\vee\>=1$ then \cite{E15} is automatically satisfied. This
follows from \cite{E14}. Most root systems have such a character
$\chi$ (see \cite{typeB} for details).

\medskip

\noindent For a $W$\_scheme $Z$ let $\|Mor|^\Phi(Z,A)$ be the set of
all $\Phi$\_equivariant morphisms $\phi:Z\pfeil A$. Again it was shown
in \cite{ARC} that the functor $Y/\fq\mapsto\|Mor|^\Phi(Y_{\fa^*},A)$
is represented by a affine commutative group scheme $\fA_\Phi/\fq$,
i.e., with
$$2
\|Mor|^\Phi(Y_{\fa^*},A)=\|Mor|_\fq(Y,\fA_\Phi).
$$

We proceed by summarizing some facts about $\fA=\fA_\Phi$. Proofs can
be found either here or in \cite{ARC}.

\medskip

\hang\noindent {\it Relation to $\fA_W$:} The group scheme $\fA_\Phi$
is an affine open subgroup scheme of $\fA_W$. In fact, all affine open
subgroup schemes are of this form. In particular, if $\Phi_{\rm max}$
is the root system where all roots are as long as possible, or
equivalently where all coroots are primitive, then $\fA_{\Phi_{\rm
max}}=\fA_W$.

\medskip

\hang\noindent {\it Global sections:} Let
$$
A^\Phi:=\{t\in A\mid\alpha(t)=1\hbox{ for all }\alpha\in\Delta\}.
$$
This is a subgroup of finite index of $A^W$, the group of $W$\_fixed points in $A$. Every $t\in A^\Phi$ induces
the constant $\Phi$\_morphism $\fa^*\pfeil A:\xi\mapsto t$ and
therefore a section $\fq\pfeil\fA$. In fact, since every morphism
$\fa^*\pfeil A$ is necessarily constant we have
$$
\Gamma(\fq,\fA)=A^\Phi.
$$

\medskip

\hang\noindent {\it Change of lattice:} Let $\Lambda'\subseteq\Lambda$ be
sublattice of finite index containing $\Delta$. Then
$\Phi'=(\Lambda',\Delta,\Lambda^{\prime\vee},\Delta^\vee)$ is also a
root datum and $A'=\|Hom|(\Lambda',\CC^\Times)$ is a quotient of $A$ by
the finite group $E:=\|Hom|(\Lambda/\Lambda',\CC^\Times)$. Because of
$\Delta\subseteq\Lambda'$ we have $E\subseteq A^\Phi$. Thus $E$ acts (freely)
on $\fA_\Phi$ by translation and $\fA_\Phi/E=\fA_{\Phi'}$. In
particular, the morphism $\fA_\Phi\pfeil\fA_{\Phi'}$ is finite and
surjective.

\medskip

\hang\noindent {\it Products of root system:} Let $\Phi$ be the
product of the root systems $\Phi_1$ and $\Phi_2$. Then
$\fa^*=\fa_1^*\times\fa_2^*$ and
$\fA_\Phi=\fA_{\Phi_1}\times\fA_{\Phi_2}$. For general root systems
one can combine this with a
suitable change of lattices as above. Thus there is always a finite
quotient $\pi:\fA\auf\fA'$ or a covering $\pi:\fA'\auf\fA$ (with $\pi$ finite)
where $\fA'$ decomposes as
$$32
\fA'=(A_0\times\fa^*_0)\times\fA_{\Phi_1}\times\ldots\times\fA_{\Phi_s}.
$$
Here, $A_0$ is a torus with Lie algebra $\fa_0$ and the $\Phi_\nu$ are
irreducible root systems.

\medskip

\hang\noindent {\it Fibers:} All fibers of $\fA/\fq$ are
commutative of dimension $r=\|dim|A=\|rk|\Phi$. More precisely, for
$a\in\fa^*$ let $s=\pi(a)$ be its image in $\fq$. Put
$$
\Delta_a:=\{\alpha\in\Delta\mid\<a\mid \alpha^\vee\>=0\}
$$
and $\Phi_a=(\Lambda,\Delta_a,\Lambda^\vee,\Delta_a^\vee)$.  Let
$\fA_s=\fA_\Phi\times_\fq\{s\}$ be the fiber over $s$ and let
$\fA_s=\fA_s^{\rm ss}\times\fA_s^{\rm uni}$ be its decomposition into
semisimple and unipotent part. Then
$$3
\fA_s^{\rm ss}\cong A^{\Phi_a}\subseteq A
\quad{\rm and}\quad
\fA_s^{\rm uni}\cong\CC^{\|rk|\<\Delta_a\>}.
$$
The first isomorphism depends on the choice of $a\in\pi^{-1}(s)$.

\medskip

\hang\noindent
{\it Generic structure:} A point $a\in\fa^*$ is called {\it regular}
if its isotropy group $W_a$ is trivial. This condition is equivalent
to $\Delta_a=\leer$, i.e., $\<a\mid \alpha^\vee\>\ne0$ for all
$\alpha\in\Delta$. Let $\fa^*_{\rm reg}$ be the set of regular points
and let $\fq_{\rm reg}$ be its image in $\fq$. Then by the description
of the fibers above we know that the restriction of $\fA$ to $\fq_{\rm
reg}$ is a torus bundle with typical fiber $A$. Its
monodromy is $W$. More precisely:
$$1
\fA|_{\fq_{\rm reg}}=A\times^W\fa^*_{\rm reg}.
$$

\medskip

\hang\noindent
{\it Local structure:} The isomorphism \cite{E1} can be generalized to
describe the local structure of $\fA$ in the neighborhood of any point
$s\in\fq$. Choose $a\in\fa^*$ lying over $s$ and let
$$
\fa^*_a:=\{b\in\fa^*\mid\<b\mid\alpha^\vee\>\ne0
\hbox{ for all }\alpha\in\Delta\setminus\Delta_a\}.
$$
This is also the set of all $b\in\fa^*$ with $\Delta_b\subseteq\Delta_a$ and
also the complement of all reflection hyperplanes which do not contain
$a$. The isotropy group $W_a$ is the Weyl group of
$\Phi_a$. Moreover, the morphism $\fq_a:=\fa^*_a/W_a\pfeil\fq$
is \'etale. Then
$$31
\matrix{\fA_{\Phi_a}|_{\fq_a}&\pfeil&\fA_\Phi\cr
\downarrow&&\downarrow\cr
\fq_a&\pfeil&\fq\cr}
$$
is a pull\_back diagram. In other words, $\fA_{\Phi_a}$ and
$\fA_\Phi$ are isomorphic over $\fq_a$. Proof: Let
$\tilde\fa^*_a:=\fq_a\times_\fq\fa^*$. Then $\fa^*_a$ is a connected
component of $\tilde\fa^*_a$. Its stabilizer in $W$ is $W(a)$ and all other
components are $W$\_translates. Now let $Y$ be any $S_a$\_scheme.
One checks
$$
\|Mor|^\Phi(Y\times_\fq\fa^*,A)=\|Mor|^\Phi(Y\times_{\fq_a}\tilde\fa^*_a,A)=
\|Mor|^{\Phi_a}(Y\times_{\fq_a}\fa^*_a,A)=\|Mor|^{\Phi_a}(Y\times_\fq\fa^*,A).
$$
Thus, the universal property implies
$\fA_\Phi\times_\fq\fq_a=\fA_{\Phi_a}|_{\fq_a}$.

\medskip

\hang\noindent {\it Lie algebra:} The Lie algebra of $\fA$ is a vector
bundle over $\fq$. Taking for $Y$ the spectrum of $\CC[X]/(X^2)$ one
obtains
$$16
\|Lie|\fA_\Phi=T^*_\fq
$$
where $T^*_{\fq}$ is the cotangent bundle of $\fq$.

\medskip

\hang\noindent
{\it Symplectic structure:} The isomorphism \cite{E16} implies that
the tangent space of $\fA$ in a point of the identity section is of
the form $T\oplus T^*$. Thus it carries a symplectic structure. We
claim, that in fact $\fA$ has a global symplectic structure, i.e.,
there is a natural closed and non\_degenerate 2-form $\omega$ on
$\fA$. It is constructed as follows. First, $A\times\fa^*$ is the
cotangent bundle of $A$, hence symplectic. The form is invariant under
the natural $W$\_action. Thus it pushes down to a symplectic structure
on $A\times^W\fa^*_{\rm reg}$. By \cite{E1}, this space is isomorphic
to the open subset $\fA_{\rm reg}:=\fA|_{\fq_{\rm reg}}$ of $\fA$ . We
show that the symplectic form on $\fA_{\rm reg}$ extends to a
non\_degenerate symplectic form on $\fA$.

\hang\hangafter0
Since $\fA$ is affine, we have to check this only in codimension one. The
complement of $\fA_{\rm reg}$ in $\fA$ is the part of $\fA$ sitting
over the ramification divisor in $\fq$. The components of this divisor
correspond to the conjugacy classes of roots. Choosing one component
and one root belonging to this component and using
\cite{E31} we are reduced to the case where the root system is of type
$\sA_1$. Using \cite{E32}, we may additionally assume that $\fA$ splits
as the product of a torus and $\fA_{\Phi'}$ where $\Phi'$ is of rank
one. Thus, we may assume that $\Phi=\Phi'$ is of rank one.

\hang\hangafter0
Then we are in the case of the example in the beginning of this
section. Recall that the coordinate ring of $\fA$ was generated by $a,b$,
and $s$ with the relation $a^2-sb^2=1$. Moreover $s=t^2$ where $t$ is the
coordinate of $\fa^*=\CC$. Let $x$ be the coordinate of
$A\cong\CC^\Times$. Then $x=a+tb$ is an eigenvalue of the
matrix \cite{E33}, i.e., the map
$$
\pmatrix{a&sb\cr b&a\cr}\mapsto(a+tb,t)\in\CC^\Times\times\CC
$$
trivializes $\fA$ over $\CC\setminus\{0\}$. A short computation yields
$$
\omega={dx\over x}\wedge dt=(a-tb)(da+b\,dt+t\,db)\wedge dt=
(a\,da-sb\,db)dt+(a\,db-b\,da)\wedge t\,dt
$$
From $a^2-sb^2=1$ it follows that the first summand vanishes. Thus
$$
\omega={1\over2}(a\,db-b\,da)\wedge ds
$$
is regular on all of $\fA$. More precisely, $\omega$ is the
restriction of a 2-form $\tilde\omega$ on $\CC^3$ to
$\fA$. Since $\fA$ is the zero\_set of $f:=a^2-sb^2-1$ and because of
$$
\tilde\omega\wedge df=(f+1)\,da\wedge db\wedge ds
$$
the restriction $\omega$ is non\_degenerate.

\beginsection ExtensionProperty. The extension property

In this section, we prove a crucial property of the group scheme
$\fA_\Phi$. It says that certain rational sections extend
automatically to regular sections. This phenomenon is related to the
group schemes introduced by Bruhat\_Tits in \cite{BrTits}~4.4. The
difference is that there the base scheme is assumed to be just
one\_dimensional.

\Definition aaaa. Let $S$ be an irreducible scheme with generic point
$\eta$ and let $\fG/S$ be a group scheme acting on an $S$\_scheme
$X/S$. Let $X_\eta$ be the fiber of $X$ over $\eta$. The action of
$\fG$ on $X$ is called {\it saturated} if for any section
$u_\eta:\eta\pfeil\fG$ of $\fG$ over $\eta$ the following are equivalent:
\Item{3} $u_\eta$ extends to a section of $\fG$ over $S$.
\Item{4} The automorphism $\phi_\eta$ of $X_\eta$ induced by $u_\eta$
extends to an $S$\_automorphism $\phi$ of $X$.\Par

\noindent Clearly, only the
implication \cite{I4}$\Rightarrow$\cite{I3} is a problem. To formulate
the results we need a restriction:

\Definition. Let $S$ be an irreducible normal affine scheme with generic point
$\eta$. An $S$\_scheme $X$ is {\it good} if it is noetherian, affine,
faithfully flat over $S$ and if the generic fiber $X_\eta$ is geometrically
irreducible.

Now we show that torus actions on good schemes are saturated. More
precisely:

\Lemma torusextend. Let $S$ be an irreducible normal $\CC$\_scheme and $X$ a
good $S$\_scheme. Let $\fT/S$ be torus acting on $X/S$. Let $\eta$ be
the generic point of $S$ and assume that the action of $\fT_\eta$ on
$X_\eta$ has a finite kernel. Then the action of $\fT/S$ on $X/S$ is
saturated

\Proof: The torus $\fT$ splits over a Galois cover $\tilde S\pfeil
S$. It is easy to see that it suffices to show the assertion for
$\tilde\fT=\fT\times_S\tilde S$ acting on $\XS=X\times_S\tilde
S$. Thus, we may assume from the outset that $\fT=T\times_\CC S$ is
split where $T$ is a torus over $\CC$.

The normality of $S$ implies that a section $u_\eta$ extends to $S$ if
and only if it extends in codimension one. Thus, after changing the
base scheme $S$, we may assume that $S=\|Spec|V$ where $V$ is a
discrete valuation ring. Let $t\in V$ be a uniformizer. Then $X=\|Spec|R$
is an affine scheme where $R$ is a faithfully flat $V$\_algebra. This
means in particular that $t$ is not a zero\_divisor in $R$, i.e., the
homomorphism $R\pfeil R[t^{-1}]$ is injective. Moreover, $X_\eta$
being irreducible means that $R[t^{-1}]$ is a domain. Hence, $R$ is a
domain, as well.

Choose an isomorphism $T\cong(\CC^\Times)^r$. Then we can write
$$
u_\eta=(t^{d_1},\ldots,t^{d_r})(v_1,\ldots,v_r)=\lambda(t)\cdot v
$$
where $d_\nu\in\ZZ$ and where $v_i\in V$ is invertible. Since $v$
is regular over $S$, we may replace $u_\eta$ by
$\lambda$. Moreover, since $u_\eta=\lambda$ takes values in a
one\_dimensional subtorus we may assume $T=\CC^\Times$ and
$u_\eta=t^d$ with $d\in\ZZ$. All we have to show now is $d=0$.

Suppose $d\ne0$. The action of $T$ on $X$ corresponds to a
$\ZZ$\_grading $R=\oplus_mR_m$. Let $\phi$ be the automorphism of
\cite{aaaa}\cite{I4}. Then $\phi(f)=t^{dm}f$ for all $f\in
R_m$. Since $\phi$ is invertible we see that multiplication by $t$ is
bijective on $R_m$ for all $m\ne0$. Put $I:=tR$ and
$J:=\bigcap_{n\ge0}I^n$. Because $R$ is a domain either $J=0$
or $I=R$ (Krull, see e.g. \cite{AtMa}~10.18). Suppose $J=0$. Since
$R_m\subseteq J$ for $m\ne0$ this would mean that $A$ is acting
trivially on $X$, contradicting the assumption that the kernel of the
action is finite. Thus $I=R$. But then $t$ is invertible in $R$,
contradicting the assumption that $X$ is {\it faithfully} flat over
$S$.\qed

In the next lemma we show that actions of $\fA_\Phi$ are almost
saturated, even in a formal sense. For this, we use two ideas: 1. the
saturatedness of tori and 2. the following simple fact: let $f(x)$ be
a rational function with $f(-x)=f(x)^{-1}$. Then $f$ is defined and
non-zero in $x=0$.

\Lemma formalextend. Let $X$ be a good scheme over
$\fq=\fa^*/W$. Assume that $\fA_\Phi$ is acting on $X$ such that the
action over the generic point of $\fq$ is faithful. Let $\hat\fq$ be
the formal neighborhood of $0$ in $\fq$ and let $\eta$ be its generic
point. Let $u$ be a section of $\fA_\Phi$ over $\eta$ such that the
action of $u$ on $X\times_\fq\eta$ extends to an automorphism of
$X\times_\fq\hat\fq$. Then $u$ extends to a section of $\fA_W$ over
$\hat\fq$.

\Proof: First observe that $\fA_\Phi$ is an open subset of $\fA_W$ and
that both coincide over $\fq_{\rm reg}$. We have to show that the
morphism $u:\eta\pfeil\fA_\Phi$ extends to a morphism
$\hat\fq\pfeil\fA_W$. Again, since $\hat\fq$ is normal, it suffices to
check this in codimension one.

Let $\hat\fq_\reg=\hat\fq\cap\fq_\reg$ be the complement of the
ramification divisor $D$. Since $\fA_\Phi$ is a torus over $\hat\fq_\reg$
we can invoke \cite{torusextend} and conclude that $u$ is regular over
$\hat\fq_\reg$.

It remains to study regularity in $D$. For a root $\alpha\in\Delta$
let $D_\alpha$ be the image of the reflection hyperplane $H_\alpha$ in
$\fq$. Then the $D_\alpha$ are precisely the irreducible components of
$D$ with $D_\alpha=D_\beta$ if and only if $W\alpha=W\beta$.

Let $\hat\fa^*$ be the formal neighborhood of $\fa^*$ in
$0$. Moreover, let $V$ be the local ring of $\hat\fa^*$ in
$\hat\fa^*\cap H_\alpha$ and let $\hat\fa^*_\alpha$ be its
spectrum. Then $V$ is a discrete valuation ring with valuation denoted
by $v_\alpha$ and uniformizer $\alpha^\vee$.

Put $W_\alpha:=\{1,s_\alpha\}$ and
$\hat\fq_\alpha:=\hat\fa^*_\alpha/W_\alpha$. Then there is a flat
morphism $\hat\fq_\alpha\pfeil\hat\fq$ with the closed point of
$\hat\fq_\alpha$ mapping to the generic point of $D_\alpha$. Thus we
have to show that the rational section $u$ is regular over
$\hat\fq_\alpha$.

\def\ahh{{\tilde u}}
\def\fahh{\hat\fa_\alpha^*}

Observe that
$$
\hat\fq_\alpha\times_\fq\fa^*=W\times^{W_\alpha}\hat\fa^*_\alpha.
$$
The section $u$ corresponds therefore to a $W_\alpha$\_equivariant
rational morphism $\xymatrix@=10pt{\tilde
u:\hat\fa^*_\alpha\ar@{.>}[r]&A}$. Now put $f:=\alpha\circ\ahh$ which
is a rational function on $\hat\fa^*_\alpha$. The
$W_\alpha$\_equivariance of $\ahh$ implies
$$8
f\circ s_\alpha=\alpha\circ\ahh\circ s_\alpha=\alpha\circ s_\alpha\circ\ahh =
(\alpha\circ\ahh)^{-1}=f^{-1}.
$$
Because the valuation $v_\alpha$ is
$s_\alpha$ invariant this implies $v_\alpha(f)=0$, i.e.,
$f$ is actually defined and invertible on $\hat\fa^*_\alpha$. Equation
\cite{E8} implies moreover that $f^2-1$ vanishes on
$H_\alpha$.

Now we consider the coroot $\alpha^\vee$ as a homomorphism $\CC^\Times\pfeil
A$ and define a morphism $b:\fahh\pfeil A$ by
$b=\alpha^\vee(f)$. Then the $W_\alpha$\_equivariance
of $\ahh$ can be rewritten as
$$
\ahh\circ s_\alpha=\ahh\cdot b^{-1}
$$
Thus
$$
c:=\ahh^2b^{-1}=\ahh\cdot(\ahh\circ s_\alpha)
$$
has values in $A^{W_\alpha}$. This group has at most two connected
components which implies that $c^2$ has values in the torus
$A_0:=(A^{W_\alpha})^0$. Now we use the fact that $u$ induces an
action on all of $X\times_\fq\hat\fq$. Thus also $c$ extends to an
automorphism of $X\times_\fq\hat\fq_\alpha$. Since it is part of a
torus action, we conclude (\cite{torusextend}) that $c$ is regular on
$\fahh$. Thus, also $\ahh^4$ and therefore $\ahh$ are regular on
$\fahh$.\qed

Before we go on, we need to study the difference between sections of
$\fA_W$ and sections of $\fA_\Phi$.  Recall that each element
$t\in A^W$ gives rise to a constant section, also denoted by $t$, of
$\fA_W$.

\Lemma globalRep. Let $\hat\fq$ be the formal neighborhood of a point
$s\in\fq$. Then
$$
\Gamma(\hat\fq,\fA_W)=A^W\cdot\Gamma(\hat\fq,\fA_\Phi).
$$

\Proof: For an algebraic group $E$ let $\pi_0(E)=E/E^0$ be
its group of components. It suffices to show that
$$24
A^W\pfeil\pi_0(\fA_{W,s})
$$
is surjective. Indeed, let $u$ be a section of $\fA_W$
over $\hat\fq$. Let $u(s)\in\fA_{W,s}$ be its value at $s$. Then there
is $t_0\in A^W$ with $u(s)\cdot
t_0^{-1}\in\fA_{W,s}^0\subseteq\fA_{\Phi,s}$. But then $u_1:=t_0^{-1}\cdot u$ is a
section of $\fA_\Phi$.

To show surjectivity of \cite{E24}, let $a\in\fa^*$ be a lift of
$s$. Because of $\pi_0(\fA_{W,s})=\pi_0(A^{W_a})$ it suffices to show
that $A^W\pfeil\pi_0(A^{W_a})$ is surjective or, equivalently, that
the map between character groups is injective.

Let $\Delta_{\rm max}\subseteq \Lambda$ be the root system with Weyl
group $W$ such that all coroots are primitive. Then $W_a$ corresponds
to a subroot system $\Delta_{\rm max}'$. Choose a system of simple
roots $\Sigma\subseteq\Delta_{\rm max}$ such that
$\Sigma'=\Sigma\cap\Delta_{\rm max}'$ is a set of simple roots for
$\Delta_{\rm max}'$. Then the map of character groups is
$$
\left(\Lambda/\<\Delta'_{\|max|}\>\right)_{\|torsion|}\Pfeil\
\Lambda/\<\Delta_{\|max|}\>.
$$
The kernel is a torsion subgroup of the free group
$\<\Delta_{\|max|}\>/\<\Delta'_{\|max|}\>\cong
\sum_{\alpha\in\Sigma\setminus\Sigma'}\ZZ\alpha$, hence zero.\qed

Next, we investigate when there is actually a difference between
$\fA_W$ and $\fA_\Phi$.  Clearly this only happens if there is a
root $\alpha$ with ${1\over2}\alpha^\vee\in\Lambda^\vee$. These roots will be called {\it special}.

\Lemma typeB. Let $\alpha$ be a special root. Then $\Phi$ has a direct
summand isomorphic to the root datum of $SO_{2n+1}(\CC)$, $n\ge1$, in
which $\alpha$ is a short root.

\Remark: The root datum of $SO_{2n+1}(\CC)$ is of course the root
system $\sB_n$ where $\Lambda$ is the root lattice. Observe that
$\sB_1=\sA_1$ is included.

\medskip

\Proof: We may choose a set of simple roots $\Sigma\subseteq\Delta$
containing $\alpha$. Then $\<\Lambda\mid \alpha^\vee\>=2\ZZ$ implies
$\<\beta|\alpha^\vee\>\in\{0,-2\}$ for all $\beta\in\Sigma$. The
classification of Dynkin diagrams shows that $\alpha$ is the short
simple root of a subroot system of type $\sB_n$. Let
$\alpha_1,\ldots,\alpha_n$ be the short positive roots. By assumption,
the corresponding half-coroots
${1\over2}\alpha_1^\vee,\ldots,{1\over2}\alpha_n^\vee$ lie in
$\Lambda^\vee$. Because of
$\<\alpha_i|{1\over2}\alpha_j^\vee\>=\delta_{ij}$, we see that the
lattice spanned by the $\alpha_i$ is a direct summand of $\Lambda$.\qed

As already mentioned, any root $\alpha$ gives rise to a component
$D_\alpha$ of the ramification divisor of $\fa^*\pfeil\fq$. We call
$D_\alpha$ {\it special} if $\alpha$ is special.

\Lemma specialDiv. Let $D\subseteq\fq$ be an irreducible divisor.
\Item{5} If $D$ is not a special divisor then $\fA_{\Phi,s}=\fA_{W,s}$
for a generic point $s\in D$.
\Item{6} If $D=D_\alpha$ is special then $\fA_{\Phi,s}$ is of index
$2$ in $\fA_{W,s}$ for a generic point $s\in D$.
\Item{1}If $D=D_\alpha$ is special then there is an involution $e_\alpha\in
A^W$ with $e_\alpha(s)\not\in\fA_{\Phi,s}$ if and only if $s\in D$.

\Proof: If $D$ is not part of the ramification divisor then, for $s\in
D$ generic, $\fA_{W,s}$ is a torus, hence
$\fA_{\Phi,s}=\fA_{W,s}$. Otherwise, $D=D_\alpha$ for some root
$\alpha$. If $\alpha$ is not special then $\alpha\in\Delta_{\rm
max}$. Thus $\fA_{\Phi,s}=\fA_{W,s}\cong\|ker|\alpha\times\CC$ for
$s\in D$ generic. Otherwise, $2\alpha\in\Delta_{\rm max}$ and
$\fA_{\Phi,s}$ is of index $2$ in $\fA_{W,s}$. This shows \cite{I5}
and \cite{I6}.

For \cite{I1} we may assume that $\alpha$ is a short root in a direct
summand $\Phi'$ of $\Phi$ which is of type $\sB_n$. Thus, there is a
unique homomorphism $\epsilon_\alpha:\Lambda\pfeil\{\pm1\}$ with
$$
\epsilon_\alpha(\chi)=\cases{-1&if $\chi$ is a short root of $\Phi'$\cr+1&if
$\chi$ is another root or if $\chi$ is fixed by $W$.\cr}
$$
This homomorphism is clearly $W$\_invariant and of order $2$. Thus, it
can be interpreted as an involution in $e_\alpha\in A^W$. Let $s\in\fq$
and let $a\in\fa^*$ lie above $s$. Then $e_\alpha(s)\in\fA_{\Phi,s}$
if and only if $\epsilon_\alpha(\Delta_a)=1$ if and only
$\Delta_a$ contains no short root of $\Phi'$ if and only if $s\not\in
D_\alpha$.\qed

Now we can give a criterion for an $\fA_\Phi$\_action to be saturated:

\Theorem formalSat. Assume $\fA_\Phi/\fq$ acts on a scheme $X/\fq$
which is good and assume that the action is generically effective. For
$s\in\fq$ let $\hat\fq_s$ be its formal neighborhood. Moreover, put
$\hat X_s:=X\times_\fq\hat\fq_s$ and
$\hat\fA_s=\fA_\Phi\times_\fq\hat\fq_s$. Then the following are
equivalent:
\Item{7}For every $s\in\fq$, the action of $\hat\fA_s$ on $\hat X_s$
is saturated.
\Item{8}The action of $\fA_\Phi$ on $X$ is saturated.
\Item{9}Let $t\in A^W$ with $t\not\in A^\Phi$. Then the rational
action of $t$ on $X$ does not extend to an action on all of $X$.

\Proof: \cite{I7}$\Rightarrow$\cite{I8}: Let $u$ be a rational section
of $\fA_\Phi$ such that the action of $u$ on $X$ is regular. Then $u$
is regular over $\hat\fq_s$ for every $s\in\fq$. Hence it is
regular on all of $\fq$.

\cite{I8}$\Rightarrow$\cite{I9}: Let $t\in A^W$ whose rational action
on $X$ is regular. Then $t$ is a regular section of $\fA_\Phi$, i.e.,
$t\in A^\Phi$.

\cite{I9}$\Rightarrow$\cite{I7}: Let $u$ be a rational section of
$\fA_\Phi$ over $\hat\fq_s$ which induces an automorphism of $\hat
X_s$. Then \cite{formalextend} implies that $u$ is a regular section
of $\fA_W$ over $\hat\fq_s$ and we have to show that $u$ is actually a
section of $\fA_\Phi$. \cite{globalRep} implies that there is $t_0\in
A^W$ such that $u_1=u\cdot t_0^{-1}$ is regular over
$\hat\fq_s$. Thus, we may replace $u$ by $t_0$.

So assume $u\in A^W$. Now let $D\subset\fq$ be an irreducible divisor
over which $u$ is not a section of $\fA_\Phi$. Then \cite{specialDiv}
tells us that $D=D_\alpha$ is a special divisor. If $s\not\in
D_\alpha$, then we may eliminate $D_\alpha$ by multiplying $u$ by the
special involution $e_\alpha$ of \cite{specialDiv}\cite{I1}. Thus, we
may assume $s\in D_\alpha$. Then the action of $u$ on $X$ is regular
outside the union of finitely many irreducible divisors which all pass
through $s$. On the other side, the action of $u$ extends to an action
on $\hat X_s$. These two facts imply that the action of $u$ on $X$ is
regular outside a subset of codimension $\ge2$. Since $X$ is affine,
the action is regular everywhere. Now \cite{I9} implies $u\in
A^\Phi$. A fortiori, $u$ is a regular section of $\fA_\Phi$ over
$\hat\fq_s$.\qed

Finally, we need that for given $W$ and $X$ at most one $\fA_\Phi$
can act in a saturated fashion.

\Theorem uniqueA. Let $\hat\fq$ be the formal neighborhood of $0$ in
$\fq$ and let $\hat X/\hat\fq$ be good. Let $\Phi_1$ and $\Phi_2$ be two root
systems with the same lattice and Weyl group. Assume that both 
$\fA_{\Phi_1}\times_\fq\hat\fq$ and $\fA_{\Phi_2}\times_\fq\hat\fq$
act on $X$ inducing the same action of their Lie algebra. Assume
moreover that both group scheme actions are saturated. Then
$\Phi_1=\Phi_2$.

\Proof: Let $E\subseteq A^W$ be the group of constant sections
which act on $X$. By saturatedness, we have
$A^{\Phi_1}=E=A^{\Phi_2}$. This implies $\Phi_1=\Phi_2$.\qed

\beginsection Spherical. Automorphisms of spherical varieties

In this section, we apply the results of the preceding section to the
cotangent bundle of spherical varieties. Our starting point is the
following theorem.

\Theorem algAuto. Let $X$ be a smooth affine spherical variety. Then
there is a root system $\Phi_X$ with lattice $\Lambda_X$ (the group
generated by the weight monoid) and Weyl group $W_X$ with:
\Item{30} The group scheme $\fA=\fA_{\Phi_X}$ acts on $Z:=T^*_X$ over
$\fq:=\fa^*/W_X$.
\Item{31}This action is uniquely determined by the requirement that the
following diagram commutes:
$$22
\cxymatrix{
T^*_\fq\times_\fq Z\ar[r]^-\beta_-\sim\ar[d]^\alpha&T^*_Z\ar[d]^\gamma_\sim\\
(\|Lie|\fA)\times_\fq Z\ar[r]^-\delta&T_Z}
$$
Here, $\alpha$ is induced by the isomorphism \cite{E16}, $\beta$ is
pull\_back of cotangent vectors via $Z\pfeil\fq$, $\gamma$ comes from
the symplectic structure on $Z$, and $\delta$ is the infinitesimal
action of $\fA$ on $Z$.
\Item{32}The action of $\fA$ on $Z$ is saturated.

\Proof: See \cite{ARC} Thm.~4.1 and Thm.~7.8.\qed

It turns out that the data $(\Lambda_X,W_X)$ determine to a large
extent the generic structure of $Z$. Let $L\subseteq G$ be the
centralizer of $\fa^*\subseteq\ft^*\subseteq\fg^*$ under the
coadjoint action. This is also the largest Levi subgroup of $G$
containing $T$ such that $\Lambda_X\subseteq X(L)$. Let $L_0$ be the
set of $g\in L$ with $\chi(l)=1$ for all $\chi\in\Lambda_X$. This is a
normal subgroup of $L$ such that $A:=L/L_0$ is a torus with Lie
algebra $\fa$.

\Proposition genericstructure. There is a dense open $W_X$\_stable subset
$\fa_0^*\subseteq\fa^*$ and an isomorphism
$$
\epsilon:G/L_0\times\fa_0^*\Pf\sim Z\times_\fq\fa_0^*
$$
with:
\Item{35}$\epsilon$ is $G$\_equivariant.
\Item{36}$\epsilon$ is compatible with the moment map, i.e., the diagram
$$
\cxymatrix{G/L_0\times\fa_0^*\ar[r]\ar[dr]_{m^0}&Z\ar[d]^m\\
&\fg^*}
$$
commutes where $m^0$ is the morphism $(gL_0,\xi)\mapsto g\xi$.
\Item{37}$\epsilon$ is compatible with the invariant moment map, i.e.,
the diagram
$$
\cxymatrix{G/L_0\times\fa_0^*\ar[r]\ar[d]_{p_2}&Z\ar[d]^\psi\\
\fa_0^*\ar[r]&\fa^*/W_X}
$$
commutes where $p_2$ is the projection to the second factor.
\Item{38}Let $W_X$ act on $Z\times_\fq\fa_0^*$ by acting on the second factor.
Then $\epsilon$ is $W_X$\_equivariant where $w\in W_X$ acts on
$G/L_0\times\fa_0^*$ according to the formula
$$40
{}^w(gL_0,\xi)=(gt_\wq(\xi)\wq^{-1}L_0,w\xi).
$$
Here, $\wq$ is a lift of $w$ to $N_G(\fa^*_X)$ and $t_\wq:\fa_0^*\pfeil
A$ is some morphism.
\Item{39}$\epsilon$ is compatible with the $\fA_X$\_action: because of
$\fA_X\times_\fq\fa_0^*=A\times\fa_0^*$ this is just an $A$\_action which
coincides with the right $A$\_action on $G/L_0$.

\Proof: Most assertions follow from \cite{CM}, especially \S4. The
main point is that the closure of $m(Z)$ coincides with the closure of
$G\fa^*$ inside $\fg^*$ (\cite{CM}~3.3). Thus
$\psi:G\times^Lm^{-1}(\fa^*)\pfeil Z$ is dominant and even generically
finite but $\psi$ is in general not birational and $m^{-1}(\fa^*)$ is
not irreducible. Let $B$ be the Borel subgroup of $G$ and $X_0=Bx_0$
be the open $B$ orbit in $X$. Then we replace $\fa^*$ by an open
subset $\fa_0^*$ (making $\psi$ finite) and $m^{-1}(\fa_0^*)$ by its
intersection $\Sigma$ with $T^*_{X_0}$. Moreover, we replace the
target $Z$ by $\tilde Z:=Z\times_{\ft^*/W}\fa_0^*$. Then $\Sigma$ is
irreducible and $\hat\psi:G\times^L\Sigma\pfeil\tilde Z$ is an open
embedding. The image is an irreducible component $\hat
T_X=T^*_X\times_{\fa^*/W_X}\fa_0^*$ of $\tilde Z$. Since $L_0$, the
principal isotropy group, acts trivially on $\Sigma$, the action
factors through $A=L/L_0$. Moreover, after possibly shrinking
$\fa_0^*$, the action of $A$ on $\Sigma$ is free and
$\Sigma\pfeil\fa_0^*$ is a quotient. Thus, again after shrinking
$\fa_0^*$, there is an isomorphism $\Sigma\cong A\times\fa_0^*$. From
this we get a map
$$
\epsilon:G/L_0\times\fa_0^*=G\times^L\Sigma\into\hat T_X\pfeil T^*_X.
$$
Now all assertion are clear by construction except for possibly
\cite{I38}.

The $W_X$\_action is of course the one on $\hat T_X$. Thus, the
$W_X$\_action commutes with the $G$\_action and the second projection
to $\fa_0^*$ is $W_X$\_equivariant. Therefore, the action of $w\in W_X$ must
have the form
$$
{}^w(gL_0,\xi)=(g n_w(\xi)L_0,w\xi)
$$
where $n_w$ is a morphism $\fa_0^*\pfeil N_G(L_0)/L_0$. Now, the fact
that $m$ is $W_X$\_invariant implies $n_w(\xi)w\xi=\xi$ for all
$\xi\in\fa_0^*$. Thus, $t_\wq(\xi):=n_w(\xi)\wq\in L/L_0=A$ as claimed.
\qed

\noindent Later on, we also need the following complements to
\cite{algAuto}:

\Corollary. With the notation above:
\Item{33} The moment map is $\fA$\_invariant, i.e., $m(az)=m(z)$ for all
$(a,z)\in\fA\times_\fq Z$.
\Item{34} The $\fA$\_action is compatible with the symplectic structure
of $Z$, i.e., the following formula holds:
$$23
\mu^*\omega_Z=pr_1^*\omega_\fA+pr_2^*\omega_Z
$$
Here $pr_1:\fA\times_\fq Z\pfeil\fA$ is the first projection and
$\mu,pr_2:\fA\times_\fq Z\pfeil Z$ the multiplication and the second
projection, respectively. Moreover, $\omega_\fA$ and $\omega_Z$ are
the symplectic forms on $\fA$ and $Z$, respectively.

\Proof: \cite{I33} follows from \cite{I36} and \cite{I39} above.

\cite{I34} As in the proof above we write
$G/L_0\times\fa_0^*=G\times^L\Sigma$ with $\Sigma\cong
A\times\fa_0^*$. Then $\Sigma$ is a Hamiltonian $A$\_variety with the
second projection as moment map. Now choose coordinates $t_\nu$ of $A$
and compatible coordinates $x_\nu$ of $\fa_0^*$. Then the symplectic
form on $\Sigma$ has the form
$$
\omega_\Sigma=\sum_\nu {dt_\nu\over t_\nu}\wedge dx_\nu+\omega_0
$$
where $\omega_0$ is a 2-form involving only the coordinates
$x_\nu$. The pull\_back $\tilde\fA$ of $\fA_X$ to $\fa_0^*$ is the
trivial group scheme with fiber $A$. Thus we have
$$
\tilde\fA\times_{\fa_0^*}\Sigma=A\times A\times\fa_0^*.
$$
Denote the coordinates of the three factors by $s_\nu$, $t_\nu$, and
$x_\nu$, respectively. Then $pr_1$, $pr_2$, and $\mu$ are, respectively,
$$
(s,t,x)\mapsto (s,x), (t,x), (st,x).
$$
A very short calculation shows that \cite{E23} holds for the
$\tilde\fA$\_action on $\Sigma$. Thus, it also holds for the
$\tilde\fA$\_action on $G\times^L\Sigma$ by $G$\_equivariance and
therefore for the $\fA_X$\_action on $Z$.\qed

Now we are able to determine the automorphisms of formal neighborhoods of closed orbits in $T^*_X$.

\Theorem formAut. Let $X$ be a smooth affine spherical variety and
$Z=T^*_X$. Let $\hat\fq$ be the formal neighborhood of a point $s$ in
$\fq$. Put $\hat Z:=Z\times_\fq\hat\fq$ and let $\phi$ be a
$G$\_equivariant automorphism of $\hat Z$ such that both the moment
map $\hat Z\pfeil\fg^*$ and the invariant moment map $\hat
Z\pfeil\hat\fq$ are $\phi$\_invariant. Then $\phi$ is induced by a
(unique) section $u$ of $\fA_X$ over $\hat\fq$.

\Proof: Let $\dot\fq$ be the generic point of $\hat\fq$ and $\dot
Z:=\hat Z\otimes_{\hat\fq}\dot\fq=Z\times_\fq\dot\fq$. We claim that
it suffices to construct a section $\dot u$ over $\dot\fq$. In fact,
the morphism $Z\pfeil\fq$ satisfies the conditions of \cite{formalSat}
(the generic fiber is irreducible by construction; see \cite{WuM}~6.6
for the faithful flatness). Thus, since the action of $\fA_X$ on $Z$
is saturated (\cite{algAuto}) the same holds for the action on $\hat
Z$ (\cite{formalSat}). Thus, $\dot u$ extends to a section $u$ over
$\hat\fq$. The automorphism induced by $u$ coincides with $\phi$
over $\dot\fq$. So it coincides over all of $\hat\fq$ since $\hat
Z\pfeil\hat\fq$ is faithfully flat.

To construct $\dot u$, let $\dot\fa^*=\dot\fq\times_\fq\fa^*$.  Thus,
via the morphism $\dot\fa^*\pfeil\fa_0^*$, the isomorphism $\epsilon$
in \cite{genericstructure} yields an isomorphism
$$
G/L_0\times\dot\fa^*\pf\sim\dot Z\times_{\dot\fq}\dot\fa^*
$$
Since $\dot Z\pfeil\dot\fq$ is $\phi$\_invariant, $\phi$ extends to an
automorphism $\tilde\phi$ of the right hand side and therefore of
$G/L_0\times\dot\fa^*$. Because $\tilde\phi$ is $G$\_equivariant and because
the projection to the second factor is $\tilde\phi$\_invariant, it follows
that there is a morphism
$t:\dot\fa^*\pfeil\|Aut|^G(G/L_0)=N_G(L_0)/L_0$ such that
$$
\tilde\phi(gL_0,\xi)=(g\,t(\xi)L_0,\xi).
$$
Now the invariance of the moment map $\dot Z\pfeil\fg^*$ implies that
$t(\xi)\xi=\xi$. Thus $t$ takes values in
$C_G(\fa^*)/L_0=L/L_0=A$. Finally, we use that $\tilde\phi$ commutes
with the $W_X$ action. From
$$0
\eqalignno{
&w\tilde\phi(eL_0,\xi)=&w(t(\xi)L_0,\xi)=
(t(\xi)\wq^{-1}t_\wq(\xi),w\xi)
\cr
&\tilde\phi w(eL_0,\xi)=&\tilde\phi(\wq^{-1}t_\wq(\xi)L_0,w\xi)=
(\wq^{-1}t_\wq(\xi)t(w\xi)L_0,w\xi)
\cr}
$$
we obtain $t(w\xi)=w\,t(\xi) w^{-1}$, i.e., $t:\dot\fa^*\pfeil A$ is
$W_X$\_equivariant. Thus $t$ descends to a section $\dot u$ of
$\fA_X$ over $\dot\fq$ inducing $\phi$.\qed

\beginsection Comparison. Comparison results

In this section, we provide some tools to pass from the algebraic to
the smooth category.

Let $K$ be a connected compact Lie group. If $V$ is a, possibly
infinite dimensional, topological representation of $K$ let
$$
V\fin:=\{v\in V\mid \|dim|\<Kv\><\infty\}
$$
be its subrepresentation of $K$\_finite vectors. It a consequence of
the Peter\_Weyl theorem that $V\fin$ is dense in $V$ whenever $V$ is locally
convex. For the next lemma observe that every $K$\_equivariant vector
bundle $\cE$ on a homogeneous space has a real algebraic
structure. Let $\Gamma^{\rm alg}$ (resp. $\Gamma^\infty$) denote
algebraic (resp. smooth) sections.

\Lemma PWvb. Let $\cE\pfeil K/H_\RR$ be an $K$\_equivariant vector
bundle of finite rank. Then there is an isomorphism
$$
\Gamma^{\rm alg}(K/H_\RR,\cE)\Pf\sim\Gamma^\infty(K/H_\RR,\cE)\fin.
$$

\Proof: Let $E$ be the fiber of $\cE$ over the base point $eH_\RR\in
K/H_\RR$. This is a finite dimensional representation of $H_\RR$ and we have
$\cE=K\times^{H_\RR}E$. Then,
$$11
\Gamma(K/H_\RR,\cE)=\|Mor|^{H_\RR}(K,E)
$$
both in the algebraic and in the smooth category. Let $U$ be any
irreducible representation of $K$. Then, by Frobenius reciprocity, a
$K$\_homomorphism $U\pfeil\Gamma(K/H_\RR,\cE)$ is the ``same'' as an
$H_\RR$\_homomorphism $\|res|_{H_\RR}U\pfeil E$. The latter is independent of
the category. Thus we have
$$
\|Hom|^K(U,\Gamma^{\rm
alg}(K/H_\RR,\cE))\Pf\sim\|Hom|^K(U,\Gamma^\infty(K/H_\RR,\cE))
$$
implying \cite{E11}.\qed

Here is now the comparison statement:

\Lemma Compare. Let $Y=Kx$ be an $K$\_orbit in a smooth real\_algebraic
$K$\_manifold $\M$. Let $\cI_Y\subseteq\cC^\infty(\M)$ be the ideal of
smooth functions vanishing in $Y$ and let $\cC^\infty_Y(\M)$ be the
$\cI_Y$\_adic completion of $\cC^\infty(\M)$ (i.e., the ring of formal power
series along $Y$). On the other hand, let $\RR[\M]$ be the ring of regular
functions on $\M$, let $I_Y=\cI_Y\cap\RR[\M]$ be the ideal of regular
functions vanishing in $Y$ and let $\RR[\![\M]\!]$ be the $I_Y$\_adic
completion of $\RR[\M]$. Then there are natural isomorphisms
$$13
\cC^\infty_Y(\M)\fin\cong
\RR[\![\M]\!]\fin\cong
\RR[\M]\otimes_{\RR[\M]^K}\RR[\![\M]\!]^K.
$$

\Proof: The second isomorphism is a theorem of Magid~\cite{Magid}.

For the first isomorphism, put $H_\RR=K_x$ and choose an $H_\RR$\_stable
complement $N$ to $\fk x$ in $T_{\M,x}$. Then Luna's slice theorem
asserts that, in the \'etale topology, a neighborhood of $Y$ in $\M$
is isomorphic to a neighborhood of $K/H_\RR$ in the fiber bundle
$K\times^{H_\RR}N$. A fortiori, this holds in the Hausdorff topology,
too. Thus, we may replace $\M$ by $K\times^{H_\RR}N$.

Let $E_n:=\oplus_{i=0}^{n-1}S^iN^*$ and $\cE_n:=K\times^{H_\RR}E_n$. Then
$\RR[\M]/I_Y^n$ and $\cC^\infty(\M)/\cI_Y^n$ are the algebraic,
resp. smooth, sections of the vector bundle $\cE_n$. Thus \cite{PWvb}
implies
$$12
\RR[\M]/I_Y^n\cong(\cC^\infty(\M)/\cI_Y^n)\fin.
$$

Now observe that for any sequence $V_n$ of $K$\_representations holds
$$
(\prod_nV_n\fin)\fin=(\prod_nV_n)\fin.
$$
This in turn implies
$$
\big(\plim{}V_n\fin\big)\fin=\big(\plim{}V_n\big)\fin
$$
and therefore
$$
\big(\plim{}\RR[\M]/I_Y^n\big)\fin=
\big(\plim{}(\cC^\infty(\M)/\cI_Y^n)\fin\big)\fin=
\big(\plim{}\cC^\infty(\M)/\cI_Y^n\big)\fin.
$$\qed

Our second comparison statement is concerned with lifting of smooth
maps. More precisely, consider the following diagram in the category
of manifolds:
$$
\cxymatrix{
&X\ar[d]^\pi\\
U\ar[r]_\alpha\ar@{.>}[ru]^\beta&Y}
$$
Given $\alpha$ and $\pi$, we say {\it $\alpha$ lifts to $X$\/} if
there is $\beta$ with $\alpha=\pi\circ\beta$. For every $u\in U$,
$x\in X$ with $y:=\alpha(u)=\pi(y)$ we obtain homomorphisms
$\hat\alpha$ and $\hat\pi$ between completions:
$$
\cxymatrix{
&\cC^\infty_x(X)\ar@{.>}[ld]_{\hat\beta}\\
\cC^\infty_u(U)&\cC^\infty_y(Y)\ar[u]_{\hat\pi}\ar[l]_{\hat\alpha}}
$$
We say, {\it $\alpha$ lifts formally to $X$\/} if for every $u\in U$
there is $x\in X$ as above and a homomorphism $\hat\beta$ with
$\hat\alpha=\hat\beta\circ\hat\pi$. Clearly, if $\alpha$ lifts then it
lifts formally. A converse is given by the following theorem:

\Theorem FormLift. We assume
\item\ri $X$ and $Y$ are real\_algebraic manifolds and $\pi$ is a
morphism.
\item\rii There is $h\in\RR[Y]$ such that $\pi$ is a closed embedding
over $Y_0:=\{y\in Y\mid h(y)\ne0\}$.
\item\riii The function $\hq:=h\circ\alpha\in\Cinf(U)$ is non\_zero and
locally analytic.\Par\noindent
Then $\alpha$ lifts if and only if it lifts formally. Moreover, the
lift $\beta$ is unique.

\Proof: The map $\beta$ is unique since by \riii{} the zero\_set of
$\hq$ is nowhere dense and $\pi$ is injective over $Y_0$. For the
existence we have to show that
$\RR[Y]\pf{\alpha^*}\Cinf(U)$ extends to
$\RR[X]\pf{\beta^*}\Cinf(U)$.  Because $\alpha$ lifts formally, hence
set theoretically, the
image of $\alpha$ is contained in the image of $\pi$. Thus, if $I$ is
the kernel of $\RR[Y]\pf{\pi^*}\RR[X]$ then
$\alpha^*(I)=0$. Therefore, we obtain a homomorphism
$\alpha^*:R:=\RR[Y]/I\pfeil\Cinf(U)$ and we can think of $R$ as a
subring of $\RR[X]$. Because of \rii, every $f\in\RR[X]$ is of the
form $g/h^N$ where $g\in R$ and $N\in\NN$. Therefore, $\beta^*(f)$ exists
if and only if $\alpha^*(g)$ is divisible by $\hq^N$. Because $\hq$ is
locally analytic, this can be checked by looking at Taylor series
(\cite{Mal}~Thm.~1.1). But for them divisibility holds since $\alpha$
lifts formally.\qed

\beginsection Automorphisms. Automorphisms of multiplicity-free spaces

In this section, we describe the automorphism group of a multiplicity free
manifold.

We start by defining a real version of the group scheme
$\fA_\Phi$. So, let $\Phi=(\Lambda,\Delta,\Lambda^\vee,\Delta^\vee)$
be a root system. Recall that $\fa^*:=\Lambda\otimes\CC$ and
$A:=\|Hom|(\Lambda,\CC^\Times)$ have the real structures
$$
\overline{\chi\otimes z}=-\chi\otimes\zq
\quad{\rm and}\quad
\aq(\chi)=\overline{a(\chi)}^{-1}
$$
respectively. Then these structures induce a
real structure on $\fA_\Phi$. Then $\fA_\Phi(\RR)$ is a Lie group
scheme over $\fq_\RR:=(\fa^*/W)(\RR)$. In the sequel, we need only the
part lying over $\fa^*_\RR/W$. Thus, we define
$$
\fA^+_\Phi:=\fA_\Phi(\RR)\otimes_{S(\RR)}\fa^*_\RR/W.
$$

Observe that $\fA^+_\Phi$ is equipped with the structural map
$\fA^+_\Phi\pfeil\fa^*_\RR/W$ such that each fiber is an (abelian) Lie
group. More precisely, there is a unit section
$\fa^*_\RR/W\pfeil\fA^+_\Phi$ and a multiplication map
$\fA^+_\Phi\times_{\fa^*_\RR/W}\fA^+_\Phi\pfeil\fA^+_\Phi$, both smooth,
making each fiber a Lie group. We are going to call such an object a
{\it Lie group scheme over $\fa^*_\RR/W$}. Let $M$ be a manifold equipped
with a smooth map $M\pfeil\fa^*_\RR/W$. Then an {\it action of
$\fA^+_\Phi$ on $M$ over $\fa^*_\RR/W$} is given by a smooth map
$\fA^+_\Phi\times_{\fa^*_\RR/W}M\pfeil M$ satisfying the usual
identities. In particular, each fiber of $\fA^+_\Phi$ over $\fa^*_\RR/W$
(an honest Lie group) acts on the corresponding fiber of $M$ over
$\fa^*_\RR/W$ in a smooth fashion.

Next, let $\fa^0$ be an affine space with $\fa^*_\RR$ as the group of
translations. Assume that $W$ also acts on $\fa^0$ in a compatible
way. Then, using a point $a\in(\fa^0)^W$, one can identify $\fa^0$
with $\fa^*(\RR)$ and therefore transport $\fA^+_\Phi$ to a Lie group
scheme over $\fa^0/W$. This scheme is independent of the choice of $a$
since $\fA^+_\Phi$ is invariant under translations by
$(\fa^*_\RR)^W$. We call it the Lie group scheme $\fA^+_\Phi$ over
$\fa^0/W$ (as opposed to $\fa^*_\RR/W$).

\Theorem groupAction. Let $M$ be a convex multiplicity free
manifold. Then there is a unique root system $\Phi=\Phi_M$ with
lattice $\Lambda_M$ and Weyl group $W_M$ and a unique action of the
Lie group scheme $\fA^+_\Phi$ over $\fq^0:=\fa^0/W_M$ on $M$ such that
\Item{45} the following diagram commutes:
$$41
\cxymatrix{
T^*_{\fq^0}\times\limits_{\fq^0} M\ar[r]^\sim\ar[d]&T^*_M\ar[d]_\sim\\
\|Lie|\fA^+_\Phi\times\limits_{\fq^0} M\ar[r]&T_M}
$$
\Item{46}Let $Y=Kx\subseteq M$ be an orbit. Then the action of
$\fA^+_\Phi$ on $\|Spec|\cC^\infty_Y(M)\fin$ is saturated.

\Proof: First of all, the action is unique since \cite{E41} prescribes
the action of the Lie algebra of $\fA^+_\Phi$. Also the root system is
unique: by construction of $W_M$, every root hyperplane $H_\alpha$
meets the momentum image $\cP$ in some point, say $a$. Let $s$ be its
image in $\fq^0$. Then in a formal neighborhood of $s$, the length
of $\alpha$ is determined by the condition of saturatedness
(\cite{uniqueA}).

For the existence, we construct the action first locally. Let
$\M=K\times^{H_\RR}(\fh_\RR^\perp\oplus S)$ be a model manifold with
corresponding spherical variety $X=G\times^HS$. Recall that there
is a root system $\Phi=\Phi_X$ such that $\fA_\Phi$ acts on
$T^*_X$. Let $n\in N_W(\fa^*)$ be as in \cite{commute} and put
$\Phi_\M:=n(\Phi)$. Then $\fA_{\Phi_\M}$ is the same as $\fA_\Phi$ but
with structure morphism twisted by
$\fa^*/W_X\pf{n\cdot}\fa^*/W_\M$. Then $\fA_{\Phi_\M}$ is acting on
$\sM$. Moreover, \cite{E21} implies that the action is defined over
$\RR$. Thus, we get an action of $\fA^+_{\Phi_\M}$ on $\M$.

Next, let $M$ be any convex multiplicity free manifold and let
$Kx\subseteq M$ be an orbit. The local cross section theorem asserts
that an open neighborhood of $Kx$ is isomorphic to
$K\times^L\Sigma_a$. Here $\Sigma$ is an open neighborhood of
$L_\RR/H_\RR$ of a model manifold $\M$ for $L_\RR$ and the subscript
$a$ means that the moment map of $\Sigma$ is shifted by $a$.  Let
$\Phi(a)$ (with $a=\psi(x)\in\cP$) be the root system attached to $\M$
and $\fa^0=a+\fa^*_\RR$. Then the Lie group scheme $\fA^+_{\Phi(a)}$ on
$\fa^0/W_\M$ acts on $\Sigma_a$, hence on a neighborhood of
$Kx$ in $M$.

Now we claim that it suffices to show that there is a root system
$\Phi_M$ with $\Phi(a)=(\Phi_M)_a$ for all $a\in\cP$. Indeed, then the
Lie group scheme $\fA^+_{\Phi_M}$ is, near $a$, isomorphic to
$\fA^+_{\Phi(a)}$ (see \cite{E31}). Moreover, by uniqueness,
\cite{E41}, the local actions glue to a global action of
$\fA^+_{\Phi_M}$.

To prove the existence of $\Phi_M$ we first remark that the local root
systems $\Phi(a)$ have the coherence property
$$42
\Phi(b)=\Phi(a)_b
$$ for all $b\in\cP$ in a neighborhood of $a$. Indeed, this follows
from \cite{uniqueA} since the action of both $\fA_{\Phi(b)}$ and
$\fA_{\Phi(a)_b}$ on a formal neighborhood of $Ky$ (with $\psi(y)=b$)
is saturated.

Let $\Phi_{\rm max}$ be the maximal root system with Weyl group $W_M$
and let $\Sigma_{\rm max}$ be a set of simple roots with respect to
the chamber $\fa^+$. We call a simple root $\alpha\in\Sigma_{\rm max}$
{\it critical} if it is not primitive in $\Lambda_M$. Critical roots
have two important properties. First, criticality is just a property
of the triple $(\alpha,\alpha^\vee,\Lambda_M)$. In particular,
$\alpha$ is critical in $\Sigma$ if and only if it is critical in any
subset of $\Sigma$ containing $\alpha$. Second, a critical simple root
is not in the same $W_M$\_orbit of any other simple root. This follows
from the classification in \cite{typeB}. Thus we can replace any set
of critical roots by their halves and still obtain a root system.

We apply this as follows: let $H_\alpha$ be the reflection hyperplane
for a critical root $\alpha$. Since $\cP$ is convex, the intersection
$H_\alpha\cap\cP$ is convex, hence connected. Thus the simple
reflection $s_\alpha$ appears in the Weyl group of $\Phi(a)$ for all
$a\in H_\alpha\cap\cP$. The coherence property \cite{E42} implies
therefore that there is $n_\alpha\in\{1,{1\over2}\}$ such that
$n_\alpha\alpha$ is a root of $\Phi(a)$ for all $a\in
H_\alpha\cap\cP$. If $\alpha$ is not critical, we define
$n_\alpha=1$. Now let $\Phi_M$ be the root system with simple
roots $\Sigma_M:=\{n_\alpha\alpha\mid\alpha\in\Sigma\}$.

It remains to show that $\Phi(a)=(\Phi_M)_a$ for all $a\in\cP$. So, let
$\alpha$ be a simple root of $\Phi(a)$. If $\alpha$ is not critical
then $n_\alpha=1$ and $\alpha$ is a simple root of $\Phi_M$. Otherwise,
$\alpha$ is a simple root of $\Phi_M$ by construction. Thus
$\Phi(a)\subseteq(\Phi_M)_a$. On the other hand, if $\alpha$ is a
simple root of $(\Phi_M)_a$ then $s_\alpha$ is in the Weyl group of
$\Phi(a)$ (\cite{globWeyl}). Thus $\alpha$ is a positive multiple of a
simple root $\alpha'$ of $\Phi(a)$. Since $\alpha'\in(\Phi_M)_a$ this
implies $\alpha=\alpha'\in\Phi(a)$.\qed

Now we are finally in the position to determine the automorphism group
of a convex multiplicity free manifold. For the statement, the smooth
invariant moment map $M\pf\psi\cP\pfeil\cP/W_M$ is denoted by
$\psi_\infty$.

\Theorem HamAuto. Let $M$ be a convex multiplicity free manifold with
momentum image $\cP$. Let $\phi:M\pfeil M$ be a diffeomorphism. Then
the following are equivalent:
\Item{24}$\phi$ is $K$\_equivariant and the momentum map is
$\phi$\_invariant, i.e.,  $m\circ\phi=m$.
\Item{25}There is a smooth section $u:\cP/W_M\pfeil\fA^+_{\Phi_M}$
which induces $\phi$, i.e. with $\phi(x)=u(\psi_\infty(x))\,x$ for all $x\in
M$.\Par\noindent
Moreover, in case these conditions hold, the following are
equivalent:
\Item{26}The map $\phi$ is an Hamiltonian automorphism of $M$, i.e.,
$\phi$ preserves additionally the symplectic structure of $M$.
\Item{27}The image of the section $u$ is a Lagrangian submanifold of
$\fA^+_{\Phi_M}$.

\Proof: Let $u$ be a section of $\fA^+:=\fA^+_{\Phi_M}$ and let $\phi$ be the
automorphism of $M$ induced by $u$. Then $\phi$ is clearly smooth and
$K$\_equivariant. Moreover, \cite{algAuto} implies that the moment map
is $\phi$\_invariant.

Conversely, let $\phi$ be a smooth $K$\_equivariant automorphism
preserving the moment map. The (continuous) map $\psi:M\pfeil\ft^+$
factors through $m$, hence it is $\phi$\_invariant as well. It
follows, that the smooth invariant moment map
$\psi_\infty:M\pfeil\fa^0/W_M$ is also $\phi$\_invariant. In
particular, $\phi$ maps each $K$\_orbit into itself.

The Lie group scheme $\fA^+$ acts freely on the dense open subset lying
over the interior of $\cP/W_M$. This implies that the section $u$ of
$\fA^+$ is unique if it exists. In particular, existence is a local
problem over $\cP/W_M$. Let $Kx\subseteq M$ be an orbit. Since the
$\fA^+$\_action is compatible with the cross section theorem (see the
proof of \cite{groupAction}), we may assume that $m(x)=0$ and that $M$
is a convex open neighborhood of $K/H_\RR$ a model variety $\M$.

Now consider the following diagram
$$
\cxymatrix{&\fA^+\times_{\fq_\RR}\M\ar[d]^{(pr_2,\mu)}\\
M\ar[r]_-{{\rm id}\times\phi}\ar@{.>}[ur]^-\gamma&\M\times_{\fq_\RR}\M}
$$
where $\mu$ denotes the $\fA^+$\_action. Since this action is
generically free, the vertical arrow is a closed embedding outside of
some hypersurface $h=0$. More precisely, $h$ is the pull\_back of a
polynomial function $h_0$ on $\fs_\RR$. Thus, since
$\psi_\infty:M\pfeil\fs_\RR$ is real analytic also the pull\_back of
$h$ to $M$ is real analytic. Therefore, \cite{FormLift} implies that
$\gamma$ exists if and only if it exists formally. For this to check,
let $Y=Ky\subseteq M$ be a second orbit. Then $\phi$ induces an
automorphism $\hat\phi$ of
$$
\cC^\infty_Y(M)\fin=\RR[M_0]\otimes_{\RR[M_0]^K}\RR[\![M_0]\!]^K
$$
(\cite{Compare}). Moreover, the automorphism of the right hand side
preserves both the moment map and the invariant moment map. Thus,
\cite{formAut} implies that $\hat\phi$ is induced by a formal section
$\hat u$ of $\fA^+$. But this means exactly that $\gamma$ exists
formally. We conclude the existence of the smooth map $\gamma$.

Now consider the composition $\delta:M\pf\beta\fA^+\times_{\fq_\RR}
\M\pf{pr_1}\fA^+$ and the diagram
$$
\cxymatrix{M\ar[r]^\delta\ar@{>>}[d]_-{\psi_\infty}&\fA^+\\
\cP/W_M\ar@{.>}[ur]_u}
$$
Then $u$ exists since $\delta$ is $K$\_invariant and $\psi_\infty$ is the
smooth quotient map by $K$. Clearly, $u$ is the section whose
existence has been asserted.

Finally, the equivalence of \cite{I26} and \cite{I27} follows from
equation \cite{E23}.\qed

\Example: A smallest example where $\fA^+_{\Phi_M}$ has non\_connected
fibers is $M=T^*_{S_2}$, the cotangent bundle of the $2$\_sphere, with
$K=SO(3)$ acting. Here $A_\RR^{\Phi_M}=\{\pm1\}$ and $-1$ corresponds
to the automorphism induced by the antipode map. For a compact
example, one can take $M=\P^1(\CC)\times\P^1(\CC)$ with $K=SU(2)$
acting. Then $-1$ induces the flip of the two factors.

\medskip

In analogy with \cite{smoothFunctions}, we state a couple of
equivalent formulations of what a section of $\fA^+$ is.
First, recall the universal property \cite{E2} of $\fA_\Phi$. For
$Y=\fA_\Phi$ the identity map corresponds to a $\Phi$\_morphism
$\fA_\Phi\times_\fq\fa^*\pfeil A$. Because $\fa^+\pfeil\fq_\RR$ is
injective, we obtain a map
$$43
\Psi:\fA_\Phi^+\Pf{\rm homeom.}\fA^+_\Phi\times_{\fq_\RR}\fa^+\Pfeil
A_\RR.
$$
For example, for the group scheme \cite{E33} this is the map
$$
(\pmatrix{a&sb\cr b&a\cr},s)\mapsto\pmatrix{a&-\sqrt{-s}\,b\cr \sqrt{-s}\,b&a\cr}
$$
(Observe, that $s\in(\i\,\RR)^2=\RR_{\le0}$.)
The map $\Psi$ fits into the commutative diagram
$$
\cxymatrix{
\cP/W_M\ar[d]^u&\cP\ar[l]^-{\rm bij.}\ar[d]^{\tilde u}\\
\fA_\Phi^+\ar[r]_\Psi&A_\RR
}
$$

Since $W_M$ does not act on $\cP$ we need the following

\Definition. Let $\tilde u:\cP\pfeil A_\RR$ be a smooth map.
\Item{}$\tilde u$ is {\it $\Phi$\_equivariant} if for every $a\in\cP$
there is a $(W_M)_a$\_invariant open neighborhood
$U\subseteq\fa^*_\RR$ and a $\Phi_a$\_equivariant smooth map $U\pfeil A_\RR$
which coincides with $\tilde u$ on $U\cap\cP$.
\Item{} $\tilde u$ is {\it formally $\Phi$\_equivariant} if for every
$a\in\cP$ the morphism $\hat\fa^*_a\pfeil A$ is
$\Phi_a$\_equivariant. Here, $\hat\fa^*_a$ is the formal neighborhood
of $a$ in $\fa^*_\RR$.

\Proposition. The assignment $u\mapsto\tilde u$ induces a bijection
between
\Item{47}smooth sections of $\fA^+_\Phi$ over $\cP$,
\Item{48}$\Phi$\_equivariant smooth maps $\cP\pfeil A_\RR$, and
\Item{49}formally $\Phi$\_equivariant smooth maps $\cP\pfeil A_\RR$.

\Proof: Let $u:\cP/W\pfeil\fA^+$ be a smooth section. Then $\tilde u$
is a composition of smooth maps:
$$
\tilde u:\cP\Pfeil\fA^+_\Phi\times_{\fq_\RR}\fa^+\Pfeil
A_\RR
$$
and therefore smooth. It is (formally) $\Phi$\_equivariant
since the right arrow is.

Conversely, let $\tilde u:\cP\pfeil A_\RR$ be formally
$\Phi$\_equivariant. To construct $u$, consider the smooth map
$$
\fA_\Phi^+\times_{\fq_\RR}\fa^+\Pfeil A_\RR\times\fa_\RR^*.
$$
Passing to the quotient by $W$ we get a smooth map $\fA_\Phi^+\pfeil
A_\RR\times^W\fa^*_\RR$ which is an isomorphism over the regular part
of $\fa^*_\RR$. The map $\tilde u$ gives rise to a map $\uq:\cP/W\pfeil
A_\RR\times^W\fa^*_\RR$. Thus, $u$ fits into the following diagram:
$$
\cxymatrix{
&\fA_\Phi^+\ar[d]\\
\cP/W\ar[r]_-\uq\ar@{.>}[ru]^-u&A_\RR\times^W\fa^*_\RR
}
$$
The lift $u$ exists formally because of the universal property of
$\fA_\Phi$. Thus $u$ exists uniquely by \cite{FormLift}.\qed

\beginsection Cohomology. Cohomology computation

For a root datum $\Phi=(\Lambda,\Delta,\Lambda^\vee,\Delta^\vee)$ with
Weyl group $W$ and a convex subset
$\cP\subseteq\fa^*=\Lambda\otimes_\ZZ\i\RR$ with non\_empty interior
let $\cL^\Phi_\cP$ be the sheaf of Lagrangian sections of $\fA^+_\Phi$
over $\cP/W$. Then \cite{HamAuto} can be rephrased as ${\rm
Aut}(M_U)=\cL^{\Phi_M}_\cP(U)$ where $U\subseteq\cP/W_M$ is open and
$M_U=\psi_\infty^{-1}(U)$. Thus, in view of \cite{localDelzant}, the
Delzant conjecture follows from the case $i=1$ of the following
theorem.

\Theorem vanCohom. Let $\Phi$ be a root datum and $\cP\subseteq\fa^+$
a locally polyhedral convex subset with non\_empty interior. Then
$H^i(\cP/W,\cL_\cP^\Phi)=0$ for $i\ge1$.

\Proof: Let $R=\<\Delta\>\subseteq\Lambda$ be the root lattice and put
$\Lambda_0:=R\oplus\Lambda^W\subseteq\Lambda$. This sublattice is of
finite index in $\Lambda$ and gives rise to a new root datum
$\Phi_0$. Then, according to section~\cite{groupscheme}, the finite
group $E:=\|Hom|(\Lambda/\Lambda_0,\CC^\Times)$ is acting freely on
$\fA^+_\Phi$ with quotient $\fA_\Phi^+/E=\fA_{\Phi_0}^+$. Since
Lagrangian sections are mapped to Lagrangian sections one gets a
short exact sequence
$$
1\pfeil E_{\cP/W}\pfeil\cL_\cP^\Phi\pfeil\cL_\cP^{\Phi_0}\pfeil1.
$$
where $E_{\cP/W}$ is the constant sheaf with fiber $E$.  Since
$\cP$ is convex, hence contractible, we have $H^i(\cP/W,E_{\cP/W})=0$ for $i\ge1$. Thus it suffices to
show that $\cL_\cP^{\Phi_0}$ has no higher cohomology.

Replacing $\Phi$ by $\Phi_0$, we may assume from now on that the root
lattice $R$ is a direct summand of $\Lambda$. Since the root lattice
of $\Delta_a$ is a direct summand of $R$ this implies that all fibers
of $\fA^+_\Phi\pfeil\fa^*_\RR/W$ are connected (see \cite{E3}). Hence,
the exponential map
$$
\|exp|:T^*_{\fa^*_\RR/W}=\|Lie|\fA^+_\Phi\pfeil\fA^+_\Phi
$$
is surjective. Being an local isomorphism, locally every section of
$\fA^+_\Phi$ is of the form $\|exp|\omega$ where $\omega$ is a 1-form
on $\fa^*_\RR/W$. Moreover, $\|exp|\omega$ is Lagrangian if and only if
$\omega$ is closed. The convexity of $\cP$ implies that then, again
locally, $\omega=df$ for some smooth function $f$ on $\cP/W$. Thus we
have shown that there is a short exact sequence
$$
0\pfeil\cK\Pfeil\cC^\infty_{\cP/W}\Pf{\|exp|df}\cL_\cP^\Phi\pfeil1.
$$
The higher cohomology of $\cC^\infty_{\cP/W}$ vanishes. Thus it
remains to show that $H^i(\cP/W,\cK)=0$ for $i\ge2$.

The map $\cP\pfeil\cP/W$ is a homeomorphism. Therefore, the pull\_back
$\cK^+$ of $\cK$ to $\cP$ has the same cohomology. For
$U^+\subseteq\cP$ open the group $\cK^+(U^+)$ consists of smooth
functions $f^+$ on $U^+$ such that $\|exp|df^+\equiv1$. Because
$\Lambda^\vee\subseteq\fa\cong\|Lie|A$ is the kernel of
$\|exp|:\fa\pfeil A$ the function $f^+$ must satisfy
$df^+\in\Lambda^\vee$ over the interior of $\cP$. This implies that,
at least locally, $f^+$ is an affine linear function
$$4
f^+(z)=\<z\mid \gamma\>+c\qquad{\rm where}\qquad\gamma\in\Lambda^\vee,
c\in\RR.
$$
Now we have to figure out when $f^+$ represents a {\it smooth}
function on $\cP/W$. By \cite{smoothFunctions} this is the case if the
Taylor series of $f^+$ in any point $a\in\cP$ is
$W_a$\_invariant. This translates into
$$
\<\alpha\mid \gamma\>=0\quad\hbox{\rm for all
$\alpha\in\Delta_a$.}
$$
Thus, we have proved: If $U^+\subseteq\cP$ is connected then
$\cK^+(U^+)$ is the set of functions $f^+(z)=\<z\mid \gamma\>+c$ where
$c\in\RR$, $\gamma\in\Lambda^\vee$ and $\<\Delta_a\mid \gamma\>=0$ for
all $a\in U^+$.

The choice of the Weyl chamber $\fa^+$ determines a system of simple
roots $\Sigma\subseteq\Delta$. Then $\Sigma_a:=\Sigma\cap\Delta_a$ is
a set of simple roots for $\Delta_a$ for all $a\in\fa^+$. Thus every
$a\in\cP$ has an open neighborhood $U^+$ such that $\cK^+(U^+)$ is the
kernel of
$$44
\Lambda^\vee\oplus\RR\Pfeil\ZZ^{|\Sigma_a|}:
(\gamma,c)\mapsto (\<\alpha|\gamma\>)_{\alpha\in\Sigma_a}.
$$
Recall that the root lattice is a direct summand of $\Lambda$. Thus
$\Sigma_a$ is part of a basis of $\Lambda$ which means that the map
\cite{E44} is surjective. Since $\alpha\in\Sigma_a$ can also be
rephrased as $a\in H_\alpha$ where $H_\alpha$ is the reflection
hyperplane corresponding to $\alpha$ we get a short exact sequence of
sheaves
$$
0\pfeil\cK^+\pfeil\Lambda^\vee_\cP\oplus\RR_\cP\pfeil\bigoplus_{\alpha\in\Sigma}
\ZZ_{H_\alpha\cap\cP}\pfeil0.
$$
Here, $\ZZ_{H_\alpha\cap\cP}$ is the sheaf which is constant $\ZZ$ on
$H_\alpha\cap\cP$ and zero outside. Also $H_\alpha\cap\cP$ is convex
for all $\alpha$ which implies $H^i(\cP,\ZZ_{H_\alpha\cap\cP})=0$ for
all $i\ge1$. From this the desired vanishing for $\cK^+$ follows.\qed

Now we are able to prove the Delzant Conjecture:

\Theorem DelCon. A convex multiplicity free manifold is uniquely
characterized by its momentum polyhedron together with its principal
isotropy group.

\Proof: The proof is now standard. Let $M$ and $M'$ be two convex
multiplicity free manifolds having the same moment polytope $\cP$ and
the same principal isotropy group $L_0$. By \cite{localDelzant} there
is a covering $U_\nu$ of $\cP$ by convex open subsets and isomorphisms
$\phi_\nu:M_{U_\nu}\pf\sim M'_{U_\nu}$.  Put $U_{\mu\nu}:=U_\mu\cap
U_\nu$, $\tilde U_\nu:=U_\nu/W_M$, and $\tilde
U_{\mu\nu}:=U_{\mu\nu}/W_M$. Then \cite{HamAuto} implies that the
automorphism $\phi_\mu^{-1}\phi_\nu$ of $M_{U_{\mu\nu}}$ is induced by
a unique section $\alpha_{\mu\nu}$ of $\cL^{\Phi_M}_\cP$ over $\tilde
U_{\mu\nu}$. By uniqueness, these sections satisfy the cocycle
condition $\alpha_{\lambda\nu}=\alpha_{\lambda\mu}\alpha_{\mu\nu}$ on
$\tilde U_\lambda\cap \tilde U_\mu\cap \tilde U_\nu$. The covering
$\tilde U_\nu$ of $\cP/W_M$ is $\cL^{\Phi_M}_\cP$\_acyclic by
\cite{vanCohom}. Hence the vanishing of
$H^1(\cP/W_M,\cL^{\Phi_M}_\cP)$ implies that there are sections
$\beta_\nu$ of $\cL^{\Phi_M}_\cP$ over $\tilde U_\nu$ with
$\alpha_{\mu\nu}=\beta_\mu^{-1}\beta_\nu$ on $\tilde U_{\mu\nu}$. Thus
replacing the automorphisms $\phi_\nu$ by $\phi_\nu\circ\beta_\nu$ one
obtains $\phi_\mu=\phi_\nu$ over $U_{\mu\nu}$. Therefore, the
$\phi_\nu$ glue to a global isomorphism $\phi:M\pf\sim M'$.\qed

\beginsection Classification. The classification of multiplicity free
manifolds

Now we are able to complete the classification of convex multiplicity
free manifolds. For this, we establish the following notation. Let
$\Phi=(\Lambda,\Delta,\Lambda^\vee,\Delta^\vee)$ be a root
datum. After choosing a system $\Delta^+$ of positive roots, these
data define the usual triple $T\subseteq B\subseteq G$, where $G$ is a
connected complex reductive group, $B$ is a Borel subgroup, and $T$ is
a maximal torus.

\Definition. Let $\Lambda_0\subseteq\Lambda$ be a subgroup. Then a
cone $\cQ\subseteq\Lambda\otimes\RR$ is called {\it multiplicity free}
for $(\Phi,\Delta^+,\Lambda_0)$ if there is a smooth affine spherical
$G$\_variety $X$ whose weight monoid $\Xi_X$ satisfies
$$
\eqalign{
\cQ={}&\hbox{convex cone generated by}\ \Xi_X\cr
\Lambda_0={}&\hbox{abelian group generated by}\ \Xi_X\cr}
$$

\medskip

\noindent Recall also the notion of a tangent cone of a subset $\cQ$
of a real vector space in $a\in\cQ$: it is the convex cone generated by
the set $\cQ-a$. The advertised classification is:

\Theorem mfClass. Let $K$ be a compact connected Lie group with root datum
$\Phi$ and a choice of positive roots $\Delta^+$. Then there is a
bijection between isomorphism classes of convex multiplicity free
$K$\_manifolds and pairs $(\cQ,\Lambda_0)$ such that
\Item{}$\Lambda_0$ is a subgroup of $\Lambda$ and $\cQ$ is a locally
polyhedral convex subset of the Weyl chamber determined by $\Delta^+$ and
\Item{2} the tangent cone of every $a\in\cQ$ is multiplicity free for
the triple $(\Phi_a,\Delta_a^+,\Lambda_0)$.

\Proof: First, note that we replaced the moment polytope $\cP$ sitting
in $\Lambda\otimes\i\,\RR$ by
$\cQ={1\over\i}\cP\subseteq\Lambda\otimes\RR$. The only thing left to
prove is the existence of $M$ for any pair $(\cQ,\Lambda_0)$
satisfying the two conditions. By definition, there is an open convex
cover $U_\nu$ of $\cP$ and convex multiplicity free manifolds $M_\nu$
with invariants $({1\over\i}U_\nu,\Lambda_0)$. Moreover, by Delzant's conjecture
these $M_\nu$ are all isomorphic over the intersection $U_\mu\cap
U_\nu$. In particular, these manifolds yield a well\_defined system of
local root systems $\Phi(a)$, $a\in\cP$. Now observe that in the
construction of the global root system $\Phi_M$, only the two coherence
properties \cite{E5} and \cite{E42} entered. Thus, there is a root
system $\Phi$ with $\Phi(a)=\Phi_a$ for all $a\in\cP$. Let $W$ be its
Weyl group.

The rest is standard. Let $U_{\mu\nu}=U_\mu\cap U_\nu$ and
$U_{\lambda\mu\nu}=U_\lambda\cap U_\mu\cap U_\nu$. The isomorphisms
$\phi_{\mu\nu}:(M_\nu)_{U_{\mu\nu}}\pf\sim(M_\mu)_{U_{\mu\nu}}$
induce an automorphism
$\phi_{\lambda\nu}^{-1}\phi_{\lambda\mu}\phi_{\mu\nu}$ of
$(M_\nu)_{U_{\lambda\mu\nu}}$ and therefore a section
$\alpha_{\lambda\mu\nu}\in\Gamma(U_{\lambda\mu\nu}/W,\cL^\Phi_\cP)$. Using
the commutativity of $\cL^\Phi_\cP$ one verifies easily that these
sections satisfy the cocycle condition
$$
\alpha_{\kappa\lambda\mu}\alpha^{-1}_{\kappa\lambda\nu}
\alpha_{\kappa\mu\nu}\alpha^{-1}_{\lambda\mu\nu}=1
$$
on $(U_\kappa\cap U_\lambda\cap U_\mu\cap U_\nu)/W$.
The cover $U_\nu/W$ is $\cL^\Phi_\cQ$\_acyclic. Hence the vanishing of
$H^2(\cQ/W,\cL^\Phi_\cQ)$ implies the existence of sections
$\beta_{\mu\nu}\in\Gamma(U_{\mu\nu},\cL^\Phi_\cQ)$ with
$\alpha_{\lambda\mu\nu}=\beta_{\lambda\mu}\beta_{\lambda\nu}^{-1}\beta_{\mu\nu}$. Thus
replacing the isomorphisms $\phi_{\mu\nu}$ by
$\phi_{\mu\nu}\beta_{\mu\nu}^{-1}$ yields $\alpha_{\lambda\mu\nu}=1$,
i.e., $\phi_{\lambda\mu}=\phi_{\lambda\mu}\phi_{\mu\nu}$. This means
that the manifolds $M_\nu$ unambiguously glue to a space $M$ over
$\cP$.\qed

It is not necessary to check condition \cite{I2} for all
$a\in\cP$. For example, condition \cite{I2} is automatically satisfied
for interior points of $\cP$. More generally, let $\cF\subseteq\cP$ be
a face. Then $\Phi_a$, $\Delta_a^+$ and the tangent cone of $\cP$ in
$a$ are the same for all points $a$ in the relative interior $\cF^0$
of $\cF$. So it suffices to check condition \cite{I2} for one point in
$\cF^0$ where $\cF$ runs through all faces of $\cP$. Moreover,
\cite{I2} is almost by construction an open condition. Thus it
suffices to check only the minimal faces. In particular, if $\cP$ is
compact then it suffices to check \cite{I2} for the vertices alone. To
do this, one only needs to consider spherical varieties whose cone
$\cQ$ is pointed. This way, we recover Delzant's theorem
\cite{DelTor}:

\Corollary. Let $K=T_\RR$ be a torus. Then the compact multiplicity
free $T_\RR$\_manifolds $M$ with $T_\RR$ acting effectively are classified
by compact simple regular polytopes $\cP$, i.e., by those polytopes
such that the tangent cone of $\cP$ in any vertex is spanned by a
basis of $X(T_\RR)$.

\Proof: Let $X$ be a smooth affine spherical $T$\_variety with
$T$. The convex cone $\cQ$ is pointed if and only if $H=T$. Thus
$X=S=\CC^r$ is a vector space on which $T$ acts with characters
$\chi_1,\ldots,\chi_r$. Sphericality of $X$ means that the $\chi_\nu$
are linearly independent while effectivity of the action implies that
$X(T)$ is generated by the $\chi_\nu$. Thus, $\cQ$ is spanned by a
basis of $X(T)$.\qed

We can also easily recover the classification of Igl\'esias \cite{Igl}:

\Corollary. Let $M$ be a compact multiplicity free
$G=SU(2)$\_manifold.Then there are the following possibilities for
$(\cP,\Lambda)$:
\Item{50} $\cP=\{x\}\subseteq\RR_{\ge0}$ and $\Lambda=0$. The
corresponding manifold $M$ is a coadjoint orbit, i.e., either a point
($x=0$) or a $\P^1_\CC$ ($x>0$).
\Item{51}$\cP=[0,y]\subseteq\RR_{\ge0}$ with $y>0$ and $\Lambda=d\ZZ$
with $d\in\{1,2,4\}$. The corresponding manifold $M$ is either
$\P^2_\CC\cong\P(\CC^2\oplus\CC)$ ($d=1$), $\P^1_\CC\times\P^1_\CC$
($d=2$), or $\P^2_\CC\cong\P(\fs\fl(2,\CC))$ ($d=4$).
\Item{52}$\cP=[x,y]\subseteq\RR_{\ge0}$ with $0<x<y$ and $\Lambda=d\ZZ$
with $d\in\ZZ_{>0}$. The corresponding manifold $M$ is a Hirzebruch
surface of degree $d$.\Par\noindent
Here we identified the weight lattice with $\ZZ$ and $\ft^+$ with
$\RR_{\ge0}$.

\Proof: We are either in case \cite{I50} or $\cP=[x,y]$, $0\le x<y$,
is a closed interval. Assume the latter. Then $\Lambda=d\ZZ$ for some
$d\in\ZZ_{>0}$.

If $x>0$ then all root systems $\Phi_a$ are empty, i.e.,
$(\cP,\Lambda)$ is a simple regular polytope. This is always the case
and leads to \cite{I52}.

Now assume $\cP=[0,y]$. Then the local model $X$ over $a=0$ is a smooth affine
spherical $SL(2,\CC)$\_variety. These are easily classified: $V=0$ and
$H=\CC^\Times$ or $H=N_H(\CC^\Times)$, leading to the cases $d=2$ or
$d=4$ in \cite{I51}, or $H=G$ and $V=\CC^2$ corresponding to $d=1$.\qed

In general, it is possible but quite tedious to give a full list of
multiplicity free cones for any given triple
$(\Phi,\Delta^+,\Lambda_0)$. Delzant \cite{Del2} has effective done
this in the rank-2-case. In general, one can use the classification of
smooth affine spherical varieties in \cite{KVS}.

On the other hand, in the past decade we have seen significant
progress on the combinatorial structure of multiplicity free cones,
especially through Luna's approach \cite{Lun} for classifying all
spherical varieties. For example, it was these methods which enabled
Losev to prove the ``Knop Conjecture'' (\cite{KnopConj}). We refer the
reader to Losev's work \cite{Los} for more information on these
matters.

Another benefit of \cite{mfClass} is that the construction of
multiplicity free manifolds is now a purely local process. For
example, consider the natural action of a torus $T=U(1)^n$ on
$X=\CC^n$. Then $X$ is multiplicity free, with $\cP=\RR_{\ge0}^n$ and
$\Lambda=\ZZ^n$. Let $x_i$ be the coordinates of $\RR^n$, choose
$\epsilon>0$ and ``cut off'' a corner of size $\epsilon$ off
$\cP$. More precisely, let
$$
\cP_\epsilon:=\{(x_i)\in\cP\mid x_1+\ldots+x_n\ge\epsilon\}.
$$
It is easily verified that $(\cP_\epsilon,\Lambda)$ is still a simple
regular polytope, so corresponds to a Hamiltonian manifold
$X_\epsilon$ (in fact, $X_\epsilon$ is simply the blow-up of $\CC^n$
in the origin). Now, \cite{mfClass} allows to perform this process
locally: let $M$ be a multiplicity free manifold with data $\cP$ and
$\Lambda$. Let $\cF$ be the smallest face of $\ft^+$ containing $\cP$
and let $a\in\cP$ be a vertex lying in the interior of $\cF$. Then
$(\cP,\Lambda)$ is a simple regular polytope near $a$, i.e., looks
like the example above. Thus, we can cut off a small enough corner of
$\cP$ at $a$ to get a polytope $\cP_\epsilon$ which corresponds to a
manifold $M_\epsilon$. This construction is not new since it
corresponds in fact to (a special case) of Lerman's symplectic cuts
\cite{Ler}. More complicated ``surgeries'' even at non\_toroidal
vertices are imaginable, though.

\rightskip0pt plus 5pt minus 5pt

\beginrefs

\B|Abk:AtMa|Sig:AtMa|Au:Atiyah, M.; Macdonald, I.|Tit:Introduction to
commutative algebra|Reihe:-|Verlag:Addison-Wesley Publishing
Co.|Ort:Reading, Mass.-London-Don Mills, Ont.|J:1969|xxx:-||

\L|Abk:BM|Sig:BiMi|Au:Bierstone, E.; Milman, P.|Tit:Composite
differentiable functions|Zs:Ann. of Math. (2)|Bd:116|S:541--558|J:1982|xxx:-||

\Pr|Abk:Bri|Sig:Bri|Au:Brion, M.|Artikel:Sur l'image de l'application
moment|Titel:S\'eminaire d'alg\`ebre Paul Dubreil et Marie-Paule Malliavin
(Paris, 1986)|Hgr:-|Reihe:Lecture Notes in
Math.|Bd:1296|Verlag:Springer|Ort:Berlin|S:177--192|J:1987|xxx:-||

\L|Abk:BrTits|Sig:BrTi|Au:Bruhat, F.; Tits, J.|Tit:Groupes r\'eductifs sur un corps local. II. Sch\'emas en groupes. Existence d'une donn\'ee radicielle valu\'ee|Zs:Publ. Math., Inst. Hautes \'Etud. Sci.|Bd:60|S:1-194|J:1984|xxx:-||

\L|Abk:Cam|Sig:Cam|Au:Camus, R.|Tit:Vari\'et\'es sph\'eriques affines lisses|Zs:Th\`ese de doctorat (Universit\'e J. Fourier)|Bd:-|S:-|J:2001|xxx:-||

\L|Abk:DelTor|Sig:\\De|Au:Delzant, T.|Tit:Hamiltoniens p\'eriodiques et images convexes de l'application moment|Zs:Bull. Soc. Math. France|Bd:116|S:315--339|J:1988|xxx:-||

\L|Abk:Del2|Sig:\\De|Au:Delzant, T.|Tit:Classification des actions hamiltoniennes compl\`etement int\'egrables de rang deux|Zs:Ann. Global Anal. Geom.|Bd:8|S:87--112|J:1990|xxx:-||

\B|Abk:GSj|Sig:GSj|Au:Guillemin, V.; Sjamaar, R.|Tit:Convexity properties of Hamiltonian group actions|Reihe:CRM Monograph Series {\bf26}|Verlag:American Mathematical Society|Ort:Providence, RI|J:2005|xxx:-||


\L|Abk:GS2|Sig:GS|Au:Guillemin, V.; Sternberg, S.|Tit:Multiplicity\_free
spaces|Zs:J. Diff. Geom.|Bd:19|S:31--56|J:1984|xxx:-||

\L|Abk:Igl|Sig:Igl|Au:Igl\'esias, P.|Tit:Les ${\rm
SO}(3)$-vari\'et\'es symplectiques et leur classification en dimension
$4$|Zs:Bull. Soc. Math. France|Bd:119|S:371--396|J:1991|xxx:-||

\B|Abk:Kac|Sig:Kac|Au:Kac, V.|Tit:Infinite dimensional Lie algebras,
3rd ed.|Reihe:-|Verlag:Cambridge University
Press|Ort:Cambridge|J:1990|xxx:-||

\L|Abk:KaLe|Sig:KaLe|Au:Karshon, Y.; Lerman, E.|Tit:The centralizer of
invariant functions and division properties of the moment
map|Zs:Illinois J. Math.|Bd:41|S:462--487|J:1997|xxx:-||

\L|Abk:Kir|Sig:\\Kir|Au:Kirwan, F.|Tit:Convexity properties of the moment
mapping. III|Zs:Invent. Math.|Bd:77|S:547--552|J:1984|xxx:-||

\B|Abk:KirBook|Sig:\\Kir|Au:Kirwan, F.|Tit:Cohomology of quotients in
symplectic and algebraic geometry|Reihe:Mathematical Notes
{\bf31}|Verlag:Princeton University Press|Ort:Princeton,
NJ|J:1984|xxx:-||

\L|Abk:WuM|Sig:\\Kn|Au:Knop, F.|Tit:Weylgruppe und Momentabbildung|%
Zs:Invent. Math.|Bd:99|S:1-23|J:1990|xxx:-||

\L|Abk:CM|Sig:\\Kn|Au:Knop, F.|Tit:The asymptotic behavior of invariant
collective motion|Zs:Invent. Math.|Bd:116|S:309--328|J:1994|xxx:-||

\L|Abk:ARC|Sig:\\Kn|Au:Knop, F.|Tit:Automorphisms, root systems, and
compactifications of homogeneous varieties|%
Zs:J. Amer. Math. Soc.|Bd:9|S:153--174|J:1996|xxx:-||

\L|Abk:KnopVerm|Sig:\\Kn|Au:Knop, F.|Tit:Towards a classification of
multiplicity free manifolds|Zs:Handout for the conference ``Journees
Hamiltoniennes'' in Grenoble |Bd:-|S:8 pages|J:Nov 29--30, 1997|xxx:-||

\L|Abk:KnopColl|Sig:\\Kn|Au:Knop, F.|Tit:Weyl groups of Hamiltonian
manifolds, I|Zs:Preprint|Bd:-|S:33 pages|J:1997|xxx:dg-ga/9712010||

\L|Abk:KnConv|Sig:\\Kn|Au:Knop, F.|Tit:Convexity of Hamiltonian manifolds|Zs:J. Lie Theory|Bd:12|S:571--582|J:2002|xxx:math/0112144||

\L|Abk:KVS|Sig:KVS|Au:Knop, F.; Van Steirteghem, B.|Tit:Classification
of smooth affine spherical
varieties|Zs:Transform. Groups|Bd:11|S:495--516|J:2006|xxx:math/0505102||

\L|Abk:Ler|Sig:Ler|Au:Lerman, E.|Tit:Symplectic cuts|Zs:Math. Res. Lett.|Bd:2|S:247--258|J:1995|xxx:-||

\L|Abk:LMTW|Sig:LMTW|Au:Lerman, E.; Meinrenken, E.; Tolman, S.; Woodward, Ch.|Tit:Nonabelian convexity by symplectic cuts|Zs:Topology|Bd:37|S:245--259|J:1998|xxx:dg-ga/9603015||

\L|Abk:Los|Sig:Los|Au:Losev, I.|Tit:Proof of the Knop conjecture|Zs:Ann. Inst. Fourier (Grenoble)|Bd:59|S:1105--1134|J:2009|xxx:math/0612561||

\L|Abk:Lun|Sig:Lun|Au:Luna, D.|Tit:Vari\'et\'es sph\'eriques de type $A$|Zs:Publ. Math. Inst. Hautes \'Etudes Sci.|Bd:94|S:161--226|J:2001|xxx:-||

\Pr|Abk:Magid|Sig:Mag|Au:Magid, A.|Artikel:Equivariant completions and tensor
products|Titel:Group actions and invariant theory (Montreal, PQ,
1988)|Hgr:-|Reihe:CMS Conf. Proc.|Bd:10|Verlag:Amer.
Math. Soc.|Ort:Providence, RI|S:133--136|J:1989|xxx:-||

\B|Abk:Mal|Sig:Mal|Au:Malgrange, B.|Tit:Ideals of differentiable functions|%
Reihe:Tata Institute of Fundamental Research Studies in Mathematics, No. 3|%
Verlag:Oxford University Press|Ort:London|J:1967|xxx:-||

\B|Abk:Mat|Sig:Mat|Au:Matsumura, H.|Tit:Commutative ring
theory|Reihe:Cambridge Studies in Advanced Mathematics {\bf8}|Verlag:Cambridge University Press|Ort:Cambridge|J:1986|xxx:-||

\L|Abk:MiFo|Sig:MiFo|Au:Mi\v s\v cenko, A.; Fomenko, A.|Tit:A
generalized Liouville method for the integration of Hamiltonian
systems|Zs:Funkcional. Anal. i Prilo\v zen.|Bd:12|S:46--56, 96|J:1978|xxx:-||

\L|Abk:Sch|Sig:Sch|Au:Schwarz, G.|Tit:Smooth functions invariant under
the action of a compact Lie
group|Zs:Topology|Bd:14|S:63--68|J:1975|xxx:-||

\L|Abk:Sja|Sig:Sja|Au:Sjamaar, R.|Tit:Convexity properties of the moment map
re-examined|Zs:Adv. Math.|Bd:138|S:46--91|J:1998|xxx:dg-ga/9408001||

\L|Abk:WoodTrans|Sig:\\Wo|Au:Woodward, Ch|Tit:The classification of transversal multiplicity-free group actions|Zs:Ann. Global Anal. Geom.|Bd:14|S:3--42|J:1996|xxx:-||

\L|Abk:Wood|Sig:\\Wo|Au:Woodward, Ch.|Tit:Spherical varieties and existence of invariant K\"ahler structures|Zs:Duke Math. J.|Bd:93|S:345--377|J:1998|xxx:-||


\endrefs

\bye